\documentclass[11pt]{article}
\usepackage{amsmath,amssymb,amsfonts}
\usepackage[latin1]{inputenc}

\usepackage{amsfonts}
\usepackage{color}
\usepackage{epsfig}
 \usepackage{dsfont}
\usepackage{amsmath}
\usepackage{hyperref}

\usepackage{graphicx} 

\usepackage{booktabs} 
\usepackage{array} 
\usepackage{paralist} 
\usepackage{verbatim} 
\usepackage{subfig} 
\usepackage{url}
\usepackage{hyperref}
\usepackage{float}

\usepackage{epsf}
\usepackage{psfrag} 
\usepackage{epsfig,float,psfrag}
\hypersetup{pdfborder={0 0 0},colorlinks,urlcolor=blue}  

\numberwithin{equation}{section}
\newtheorem{thm}{Theorem}[section]
\newtheorem{aadef}[thm]{Definition}
\newtheorem{alem}[thm]{Lemma}
\newtheorem{aprop}[thm]{Proposition}

\newenvironment{adem}[1][]%
   {\ \\ {\bf Proof #1~: }}%
   {\hfill\mbox{\rule{2 true mm}{3 true mm}}\vskip 2 ex\noindent}
\newenvironment{_adem}[1][]%
   {\ \\ {\bf Proof #1}}%
   {\hfill\mbox{\rule{2 true mm}{3 true mm}}\vskip 2 ex\noindent}
   {\ \\ {\bf Example #1~: }}%
   {\hfill\mbox{\rule{2 true mm}{3 true mm}}\vskip 2 ex\noindent}

\title{Asymptotics for the normalized error of the Ninomiya-Victoir scheme}
\author{A. Al Gerbi, B. Jourdain\thanks{Universit\'e Paris-Est, Cermics (ENPC), INRIA, F-77455, Marne-la-Vall\'ee, France
    e-mails: jourdain@cermics.enpc.fr, anis.al-gerbi@cermics.enpc.fr - This research benefited
    from the support of the ``Chaire Risques Financiers'', Fondation du
    Risque.}~ and E. Cl\'ement\thanks{Universit\'e Paris-Est, LAMA (UMR 8050), UPEMLV, UPEC, CNRS, F-77454, Marne-la-Vall\'ee, France,
    e-mail: emmanuelle.clement@u-pem.fr.}}

\begin{document}
\maketitle

In \cite{AJC} we proved strong convergence with order $1/2$ of the Ninomiya-Victoir scheme $X^{NV,\eta}$ with time step $T/N$ to the solution $X$ of the limiting SDE. 
In this paper we check that the normalized error defined by $\sqrt{N}\left(X - X^{NV,\eta}\right)$ converges to an affine SDE with source terms involving the Lie brackets between the Brownian vector fields. The limit does not depend on the Rademacher random
variables $\eta$. This result can be seen as a first step to adapt to the Ninomiya-Victoir scheme the central limit theorem of Lindeberg Feller type, derived in \cite{BA_K} for the multilevel Monte Carlo estimator based on the Euler scheme. When the Brownian vector fields commute, the limit vanishes. This suggests that the rate of convergence  
is greater than $1/2$ in this case and we actually prove strong convergence with order $1$.

\section{Introduction}

We consider a general $n$-dimensional stochastic differential equation, driven by a $d$-dimensional standard Brownian motion $W = \left(W^1,\ldots,W^d\right)$, of the form
\begin{equation}
\left\{
    \begin{array}{ll}
dX_t = b(X_t) dt + \sum \limits_{j=1}^d \sigma^j(X_t)dW_t^j, ~	 t \in [0,T]\\
X_0 = x
\end{array}
\right.
\label{EDS_ITO}
\end{equation}
where $x \in \mathbb{R}^n$ is the starting point, $b: \mathbb{R}^n \longrightarrow \mathbb{R}^n$ is the drift coefficient and $\sigma^j: \mathbb{R}^n \longrightarrow \mathbb{R}^n, j \in \left\{1,\ldots,d\right\}$, are the Brownian vector fields.
We are interested in the study of the normalized error process for the Ninomiya-Victoir scheme. To do so we will consider in the whole paper a regular time discretization, with time step $h = T/N$, of the time interval $[0,T]$. We introduce some notations to define the Ninomiya-Victoir scheme. Let
\begin{itemize}
\item $\left(t_k = k h \right)_{k \in [\![0;N]\!]}$ be the subdivision of $[0,T]$ with equal time step h,
\item $\Delta W^j_s = W^j_s -W^j_{t_k}$, for $s \in \left(t_k,t_{k+1}\right]$ and $j \in \left\{1,\ldots,d\right\}$,
\item $\Delta s = s - t_k$, for $s \in \left(t_k,t_{k+1}\right]$,
\item $\eta = \left(\eta_k\right)_{k \ge 1}$ be a sequence of independent, identically distributed Rademacher random variables independent of $W$.
\end{itemize}
For $V: \mathbb{R}^n \longrightarrow \mathbb{R}^n$ Lipschitz continuous, $\exp(tV)x_0$ denotes the solution, at time $t \in \mathbb{R}$, of the following ordinary differential equation in $\mathbb{R}^n$  
\begin{equation}
\left\{
    \begin{array}{ll}
 \frac{dx(t)}{dt}  = V\left(x(t)\right) \\
  x(0) = x_0. 
\end{array}
\right.
\label{ODE}
\end{equation}

To deal with the Ninomiya-Victoir scheme, it is more convenient to rewrite the stochastic differential equation \eqref{EDS_ITO} in Stratonovich form. Assuming $\mathcal{C}^{1}$ regularity for the vector fields, the Stratonovich form of \eqref{EDS_ITO} is given by:
\begin{equation}
\left\{
    \begin{array}{ll}
dX_t = \sigma^0(X_t) dt + \sum \limits_{j=1}^d \sigma^j(X_t)\circ dW_t^j \\
X_0 = x
\end{array}
\right.
\label{EDS_STO}
\end{equation}
where $\sigma^0 = b - \frac{1}{2} \sum \limits_{j=1}^d \partial \sigma^j \sigma^j $ and $\partial \sigma^j$ is the Jacobian matrix of $\sigma^j$ defined as follows
\begin{equation}
\label{Jacobian}
\partial \sigma^j = \left(\left(\partial \sigma^j \right)_{ik}\right)_{i,k \in [\![1;n]\!] } = \left(\partial_{x_k} \sigma^{ij} \right)_{i,k \in [\![1;n]\!] }. 
\end{equation} 
Now, we present the Ninomiya-Victoir scheme introduced in \cite{NV}.
\begin{itemize}
\item Starting point: $X^{NV,\eta}_{t_0} = x$.
\item For $k \in \left\{0\ldots,N-1\right\}$,
 if $\eta_{k+1} = 1 $:
\begin{equation}
 X^{NV,\eta}_{t_{k+1}} = \exp\left(\frac{h}{2}\sigma^0\right) \exp\left (\Delta W^d_{t_{k+1}}\sigma^d \right) \ldots \exp\left (\Delta W^1_{t_{k+1}}\sigma^1 \right)  \exp\left(\frac{h}{2}\sigma^0\right) X^{NV,\eta}_{t_{k}}, 
\label{case 1 NV}
\end{equation}
and if $\eta_{k+1} = -1 $:
\begin{equation}
X^{NV,\eta}_{t_{k+1}} = \exp\left(\frac{h}{2}\sigma^0\right) \exp\left (\Delta W^1_{t_{k+1}}\sigma^1 \right) \ldots \exp\left (\Delta W^d_{t_{k+1}}\sigma^d \right)  \exp\left(\frac{h}{2}\sigma^0\right)  X^{NV,\eta}_{t_{k}}.
\label{case 2 NV}
\end{equation}
\end{itemize}
The strong convergence properties of a numerical scheme, which approximates the diffusion \eqref{EDS_ITO}, are useful to control the variance of the multilevel Monte Carlo estimator based on this scheme (see \cite{Giles} and \cite{PL}). This motivated our study of the  strong convergence of the Ninomiya-Victoir scheme in \cite{AJC}. More precisely, under some regularity assumptions on the coefficients of the SDE, we proved strong convergence with order $1/2$:
\begin{equation}
\forall p \ge 1, \exists C_{NV} \in \mathbb{R}_+^*, \forall N \in \mathbb{N}^*, \mathbb{E}\left[ \underset{0 \leq k\leq N}{\max}\left\|X_{t_k} - X^{NV,\eta}_{t_k}\right\|^{2p} \right] \leq C_{NV} \left(1+ \left\| x\right\|^{2p}\right) h^p.
\end{equation}
In this present paper, we focus on the convergence in law of the normalized error defined by $\sqrt{N}\left(X - X^{NV,\eta}\right)$. The asymptotic distribution of the normalized error for the continuous time Euler scheme was established by Kurtz and Protter in \cite{KP}. The asymptotic behavior of the normalized error process for the continuous time Milstein scheme \cite{Milstein}, which is known to exhibit strong convergence with order $1$, was studied by Yan in \cite{Yan}. In both cases, the normalized error converges to the solution of an affine SDE with a source term involving additional randomness given by a Brownian motion independent of the one driving both the SDE and the scheme. 

This paper is organized as follows. In Section 2, we recall basic facts about the theory of stable convergence in law, introduced by R\'enyi \cite{Renyi} and developed by Jacod \cite{Jacod} and Jacod-Protter \cite{JP}. In Section 3, we will discuss the interpolation between time grid points and then derive the asymptotic error distribution for the Ninomiya-Victoir scheme in the general case. More precisely, we prove the stable convergence in law of $\sqrt{N}\left(X - X^{NV,\eta}\right)$ to the solution of the following SDE:
\begin{equation*}
V_t = \sqrt{\frac{T}{2}} \sum \limits_{j=1}^d \sum \limits_{m=1}^{j-1}  \int_0^t \displaystyle \left[\sigma^j,\sigma^m\right]\left(X_s\right) dB^{j,m}_s + \int_0^t \displaystyle \partial b\left(X_s \right)V_s ds + \sum \limits_{j=1}^d \int_0^t \displaystyle \partial \sigma^j\left(X_s \right)V_s dW_s^j,
\end{equation*} 
where $\left[\sigma^j,\sigma^m\right] = \partial \sigma^m \sigma^j - \partial \sigma^j \sigma^m $, for $j,m \in \left\{1,\ldots,d\right\}, m < j$, denotes the Lie bracket between the Brownian vector fields $\sigma^j$ and $\sigma^m$, $\partial b$ is the Jacobian matrix of $b$, defined analogously to \eqref{Jacobian}, and $\left(B_t\right)_{0\leq t \leq T}$ is a standard $\frac{d(d-1)}{2}$-dimensional Brownian motion independent of $W$. This result ensures that the strong convergence rate is actually $1/2$. Moreover, it can be seen as a first step to adapt to the Ninomiya-Victoir scheme the central limit theorem of Lindeberg Feller type, derived by Ben Alaya and Kebaier in \cite{BA_K} for the multilevel Monte Carlo estimator based on the Euler scheme. Their approach leads to an accurate description of the optimal choice of the parameters for the multilevel Monte Carlo estimator. When the Brownian vector fields commute, the limit vanishes, which suggests that the rate of convergence  is greater than $1/2$.
In Section 4, we focus on the commutative case and we provide a suitable interpolation between time grid points, to show strong convergence with order $1$.

\section{Stable convergence}

We start with the definition of the stable convergence in law which is stronger than the convergence in law.
\begin{aadef}
Let $\left(Z^N\right)_{N \in \mathbb{N}}$ be a sequence of random variables all defined on the same probability space $\left(\Omega,\mathcal{F},\mathbb{P}\right)$ and with values in a metric space $\left(E,d\right)$.  
Let $\left(\Omega^*,\mathcal{F}^*,\mathbb{P}^*\right)$ be an "extension" of $\left(\Omega,\mathcal{F},\mathbb{P}\right)$, and let $Z$ be an $E$-valued variable on this extension. The sequence $\left(Z^N\right)_{N \in \mathbb{N}}$ stably converges in law to $Z$ and we write this convergence as follows
\begin{equation*}
Z^N  \overset{stably}{\underset{N \to +\infty}{\Longrightarrow}} Z
\end{equation*} 
if, and only if, for all $f : E \longrightarrow \mathbb{R}$ bounded continuous and for all bounded random variable $\Xi$ on $\left(\Omega,\mathcal{F},\mathbb{P}\right)$:
\begin{equation*}
\mathbb{E}\left[f\left(Z^N\right)\Xi\right] \underset{N \to +\infty}{\longrightarrow} \mathbb{E}^*\left[f\left(Z\right)\Xi\right].
\end{equation*} 
\end{aadef}
We do not go into details of the definition of an "extension" (see \cite{Jacod} for more information).
The purpose of this section is to recall basic facts about stable convergence to study a sequence of stochastic differential equations in $\mathbb{R}^n$ of the form
\begin{equation}
\label{method_AED}
U^{N}_t = R^N_t + J^N_t + \left( \int_0^t \displaystyle  H^{0,N}_s U^N_s ds + \sum \limits_{j=1}^d \int_0^t \displaystyle H^{j,N}_s U^N_s dW_s^j \right),
\end{equation} 
where $H^{j,N}, j \in  \left\{0,\ldots,d\right\}$, take values in $\mathbb{R}^n \times \mathbb{R}^n$, $R^N$ is a remainder term and $J^N$ a source term.
This is motivated by the decomposition of the error process \eqref{Decomp_}.

The following fundamental proposition will be used to study the stable convergence in law of a random sequence of couple of variables (see section 2-1 in \cite{Jacod}).
\begin{aprop}
\label{Basic}
Let $\left(\Lambda^N\right)_{N \in \mathbb{N}}$ and $\left(\Gamma^N\right)_{N \in \mathbb{N}}$  be two sequences of random variables  all defined on the same probability space $\left(\Omega,\mathcal{F},\mathbb{P}\right)$, with values in a metric space $\left(E,d\right)$, and $\Lambda$ be a random variable on an extension, with values in $\left(E,d\right)$.
Let $\left(\Theta^N\right)_{N \in \mathbb{N}}$ be a sequence of random variables and $\Theta$ be a random variable all defined on $\left(\Omega,\mathcal{F},\mathbb{P}\right)$, with values in an other metric space $\left(E^{\prime},d^{\prime}\right)$. Then
\begin{enumerate}[(i)]
\item
 \begin{equation}
 \label{i}
\text{if }  \Lambda^N \overset{stably}{\underset{N \to +\infty}{\Longrightarrow}} \Lambda \text{ and }  d\left(\Lambda^N,\Gamma^N\right) \overset{\mathbb{P}}{\underset{N \to +\infty}{\longrightarrow}} 0 \text{ then } \Gamma^N \overset{stably}{\underset{N \to +\infty}{\Longrightarrow}} \Lambda, 
\end{equation}
\item 
\begin{equation}
\label{2ii}
\text{if }  \Lambda^N \overset{stably}{\underset{N \to +\infty}{\Longrightarrow}} \Lambda \text{ and }   d^{\prime}\left(\Theta^N,\Theta\right) \overset{\mathbb{P}}{\underset{N \to +\infty}{\longrightarrow}} 0 \text{ then } \left(\Lambda^N,\Theta^N\right) \overset{stably}{\underset{N \to +\infty}{\Longrightarrow}} \left(\Lambda,\Theta\right),
 \end{equation}
for the product topology on $E\times E^{\prime}$.
\end{enumerate}
\end{aprop}

In the following, we work on the filtered probability space $\left(\Omega,\mathcal{F}, \mathbb{F}, \mathbb{P}\right)$, where $\mathbb{F} = \left(\mathcal{F}_t = \sigma\left(\eta, W_s, s \leq t\right)\right)_{t \in [0,T]}$.
We consider the metric space $E =  \mathcal{C}\left([0,T],\mathbb{R}^n\right)$ equipped with the supremum-norm.
The following theorem, dedicated to the convergence of a sequence of semimartingales, is a simplified version of Theorem 2.1 in \cite{Jacod}.
\begin{thm}
\label{Jacod}
Let $\left(Y^N\right)_{N \in \mathbb{N}}$ be a sequence of continuous semimartingales with values in $\mathbb{R}^p$, such that $Y^N = M^N + A^N, \forall N \in \mathbb{N}$, where $M^N$ is a sequence of continuous $\mathbb{F}$-local martingales null at $t=0$ and $A^N$ is a sequence of  $\mathbb{F}$-predictable continuous processes with finite variation.  
Assume that, there exist $A$ and $f$ such that:
\begin{enumerate}
\item
\begin{equation}
\underset{t\leq T}{\sup}\left\| A^N_t - A_t \right\| \overset{\mathbb{P}}{\underset{N \to +\infty}{\longrightarrow}} 0,
\end{equation}
\item
\begin{equation}
\forall i,j \in \left\{1,\ldots,p\right\}, \forall t \in [0,T], \left<M^{i,N}, M^{j,N} \right>_t  \overset{\mathbb{P}}{\underset{N \to +\infty}{\longrightarrow}} F^{ij}_t = \displaystyle\int_{0}^{t} f^{ij}_s ds,
\end{equation}
\begin{equation}
\forall i \in \left\{1,\ldots,p\right\},\forall k \in \left\{1,\ldots,d\right\}, \forall t \in [0,T], \left<M^{i,N}, W^k \right>_t  \overset{\mathbb{P}}{\underset{N \to +\infty}{\longrightarrow}} 0.
\end{equation}
\end{enumerate}
Then, 
\begin{equation}
Y^N \overset{stably}{\underset{N \to +\infty}{\Longrightarrow}} Y
\end{equation}
where 
\begin{equation}
Y_t = A_t +  \displaystyle\int_{0}^{t} \left(f_s\right)^\frac12 dB_s, 
\end{equation}
$\left(f_s\right)^{\frac12}$ is the square root of the positive semi-definite matrix $f_s  = \left(f_s^{ij}\right)_{i,j \in [\![1;p]\!] }$ and $B$ a $p$-dimensional standard Brownian motion defined on a Wiener space $\left(\Omega^B,\mathcal{F}^B,\mathbb{P}^B\right)$ and independent of $W$. The stable convergence takes place in the canonical Wiener extension of $W$, denoted by $\left(\Omega^*,\mathcal{F}^*,\mathbb{P}^*\right)$ defined as follows
\begin{equation*}
\Omega^* = \Omega \times \Omega^B, ~\mathcal{F}^* = \mathcal{F} \otimes \mathcal{F}^B, ~\mathbb{P}^* = \mathbb{P} \otimes \mathbb{P}^B.
\end{equation*}
\end{thm}
In comparison with Theorem 2.1 of \cite{Jacod}, the assumption, 
\begin{equation*}
\left<M^i, N \right>_t  = 0, \forall i \in \left\{1,\ldots,p\right\} \text{ and } N \text{ a bounded martingale orthogonal to } W,
\end{equation*}
is obvious, since we can write $M$ in terms of an It\^o integral with respect to the Brownian motion $W$, by using the martingale representation theorem. We will use Theorem \ref{Jacod}, together with the following proposition to study the source term $J^N$ in the decomposition \eqref{method_AED}. This proposition is a consequence of Theorem 2.3  in \cite{JP} (see the proof of Theorem 2.5 (c) in \cite{JP}).                                           
\begin{aprop}
\label{Jacod_Protter}
Let $\left(Y^N\right)_{N \in \mathbb{N}}$ be a sequence of continuous semimartingales with values in $\mathbb{R}^p$, such that $Y^N_t = Y^N_0 + M^N_t + A^N_t, \forall N \in \mathbb{N}, \forall t \in \left[0,T\right] $, where $M^N$ is a sequence of continuous $\mathbb{F}$-local martingales null at $t=0$ and $A^N$ is a sequence of $\mathbb{F}$-predictable continuous processes with finite variation null at $t=0$.  
Assume that the sequence $\left(\left<M^N\right>_T + \displaystyle\int_{0}^{T} \left| dA^N_s\right|\right)_{N\in \mathbb{N}}$ is tight. 
Then, for any sequence $\left(K^N\right)_{N \in \mathbb{N}}$ of $\mathbb{F}$-predictable, right-continuous and left-hand limited processes, with values in $\mathbb{R}^q \times \mathbb{R}^p$,  such that the sequence $\left(K^N, Y^N\right)$ stably converges in law to a limit $\left(K,Y\right)$ we have the following result:\\
$Y$ is a semimartingale and with respect to the filtration generated by the limit process $\left(K,Y\right)$ and
\begin{equation}
\left(K^N,Y^N, \displaystyle\int K^N dY^N\right) \overset{stably}{\underset{N \to +\infty}{\Longrightarrow}} \left(K,Y, \displaystyle\int K dY\right),
 \end{equation}
where $\displaystyle\int K^N dY^N = \left(\displaystyle\int_0^t K^N_s dY^N_s\right)_{t\in [0,T]}$ and $\displaystyle\int K dY = \left(\displaystyle\int_0^t K_s dY_s\right)_{t\in [0,T]}$.
\end{aprop}

The following theorem deals with a sequence of stochastic differential equations in $\mathbb{R}^n$ of the form
\begin{equation}
U^N_t = R^N_t + J^N_t + \sum \limits_{j=0}^d \displaystyle\int_{0}^{t} H^{j,N}_s U^N_s dW^j_s
\end{equation}
where, by convention, $dW^0_s = ds$, $\left(J^N\right)_{N\in \mathbb{N}}$ is a sequence of continuous adapted processes, and for  $j \left\{0,\ldots,d\right\}, \left(H^{j,N}\right)_{N\in \mathbb{N}}$ is a sequence of $\mathbb{F}$-predictable, right-continuous and left-hand limited processes, with values in $\mathbb{R}^n \times \mathbb{R}^n$.

\begin{thm}
\label{Jacod_Protter_EDS}
Assume that there exist $\left(H^j\right)_{0\leq j \leq d}$ and $J$ such that:
\begin{itemize}
\item $\forall j \in \left\{0,\ldots,d\right\}, \underset{t\leq T}{\sup} \left\| H^{j,N}_t - H^j_t \right\| \overset{\mathbb{P}}{\underset{N \to +\infty}{\longrightarrow}} 0$.
\item $J^N \overset{stably}{\underset{N \to +\infty}{\Longrightarrow}} J$.
\item $\underset{t\leq T}{\sup} \left\| R^N_t \right\| \overset{\mathbb{P}}{\underset{N \to +\infty}{\longrightarrow}} 0$.
\end{itemize}
Then, $U^N$ stably converges in law towards $U$, where $U$ is the unique solution of the following affine stochastic differential equation:
\begin{equation}
U_t = J_t + \sum \limits_{j=0}^d \displaystyle\int_{0}^{t} H^{j}_s U_s dW^j_s.
\end{equation}
\end{thm}
\begin{adem}
On the one hand, denoting by 
\begin{equation*}
V^N_t = \sum \limits_{j=0}^d \displaystyle\int_{0}^{t} H^{j,N}_s dW^j_s,
\end{equation*}
the first assumption ensures that
\begin{equation*}
\underset{t\leq T}{\sup} \left\| V^N_t -V_t \right\| \overset{\mathbb{P}}{\underset{N \to +\infty}{\longrightarrow}} 0
\end{equation*}
where 
\begin{equation*}
V_t = \sum \limits_{j=0}^d \displaystyle\int_{0}^{t} H^{j}_s dW^j_s.
\end{equation*}
On the other hand, \eqref{i} from Proposition \ref{Basic} gives us 
\begin{equation*}
R^N + J^N \overset{stably}{\underset{N \to +\infty}{\Longrightarrow}} J.
\end{equation*}
Then, applying \eqref{2ii} from Proposition \ref{Basic}, we have
\begin{equation*}
 \left(R^N + J^N ,V^N\right) \overset{stably}{\underset{N \to +\infty}{\Longrightarrow}} \left(J,V\right).
\end{equation*}
Finally, since $\left(\underset{t\leq T}{\sup} \left\|H^{N}_t \right\| \right)_{N \in \mathbb{N}^*}$ is tight, we get the desired result using Theorem 2.5 (c) in \cite{JP}.  
\end{adem}

\section{Asymptotic error distribution for the Ninomiya-Victoir scheme in the general case}
\subsection{Main result}
In order to study the stable convergence in law of the normalized error process, we have to build an interpolated scheme. Let us first introduce some more notation.
\begin{itemize}
\item Let $\hat{\tau}_s$ be the last time discretization before $s \in [0,T] $, ie $\hat{\tau}_s = t_k$ if $s \in \left(t_k,t_{k+1}\right]$, and for $s = t_0 = 0 $, we set $\hat{\tau}_0 = t_{0}$.
\item Let $\check{\tau}_s$ be the first time discretization after $s \in [0,T]$, ie $\check{\tau}_s = t_{k+1}$ if $s \in \left(t_k,t_{k+1}\right]$, and for $s = t_0 = 0 $, we set $\check{\tau}_0 = 0$.
\item By a slight abuse of notation, we set $\eta_s = \eta_{k+1}$ if $s \in (t_k,t_{k+1}]$, 
\end{itemize}
A natural and adapted interpolation, at time $t\in [0,T]$, for the Ninomiya-Victoir scheme could be defined as follows:
\begin{equation}
 h_{\eta_t}\left(\frac{\Delta t}{2},\Delta W_t, \frac{\Delta t}{2}; X^{NV,\eta}_{\hat{\tau}_t} \right), 
\end{equation} 
where $\Delta W_t = \left(\Delta W_t^1, \ldots , \Delta W_t^d \right)$,
\begin{equation}
 h_{-1}\left(t_0,\ldots,t_{d+1};x\right) = \exp\left(t_0 \sigma^0\right) \exp\left(t_1 \sigma^1\right) \ldots \exp\left(t_d \sigma^d\right)  \exp\left(t_{d+1} \sigma^0\right)  x,  
\end{equation} 
and 
\begin{equation}
 h_{1}\left(t_0,\ldots,t_{d+1};x\right) = \exp\left(t_0 \sigma^0\right) \exp\left(t_d \sigma^d\right) \ldots \exp\left(t_1 \sigma^1\right)  \exp\left(t_{d+1} \sigma^0\right)  x.  
\end{equation} 
Here, to compute the It\^o decomposition of $\left(h_{\eta_t}\left(\frac{\Delta t}{2},\Delta W_t, \frac{\Delta t}{2}; X^{NV,\eta}_{\hat{\tau}_t} \right)\right)_{t\in [0,T]}$ the main difficulty is to explicit the derivatives of $h_{\eta}$. In the
general case, the computation of derivatives of this function is quite complicated.
For this reason, in this paper, we use the interpolation of the Ninomiya-Victoir introduced in \cite{AJC}:
\begin{equation}
\left\{
    \begin{array}{ll}
dX^{NV,\eta}_t = \displaystyle  \sum \limits_{j=1}^d \sigma^j(\bar{X}^{j,\eta}_t)  dW_t^j + \displaystyle \frac{1}{2} \sum \limits_{j=1}^d \partial \sigma^j \sigma^j\left(\bar{X}^{j,\eta}_t\right) dt + \frac{1}{2} \left(\sigma^0\left(\bar{X}^{0,\eta}_t\right) + \sigma^0\left(\bar{X}^{d+1,\eta}_t\right)\right) dt   \\
X^{NV,\eta}_0 = x
\end{array}
\right.
\label{NV-Interpol_ITO}
\end{equation}
where, for $s \in \left(t_k, t_{k+1}\right]$:
\begin{equation}
\bar{X}^{0,\eta}_s = \exp\left(\frac{\Delta s}{2}\sigma^0\right) \left( X^{NV,\eta}_{t_{k}} \mathds{1}_{\left\{\eta_{k+1} = 1\right\}} + \bar{X}^{1,\eta}_{t_{k+1}} \mathds{1}_{\left\{\eta_{k+1} = -1\right\}} \right),
\end{equation}
\begin{equation}
\text{for }j \in \left\{1,\ldots,d\right\}, \bar{X}^{j,\eta}_s = \exp\left(\Delta W_s^j\sigma^j\right) \left( \bar{X}^{j-1,\eta}_{t_{k+1}} \mathds{1}_{\left\{\eta_{k+1} = 1\right\}} + \bar{X}^{j+1,\eta}_{t_{k+1}} \mathds{1}_{\left\{\eta_{k+1} = -1\right\}} \right),
\label{Bar_j}
\end{equation}
\begin{equation}
\bar{X}^{d+1,\eta}_s = \exp\left(\frac{\Delta s}{2}\sigma^0\right) \left( \bar{X}^{d,\eta}_{t_{k+1}} \mathds{1}_{\left\{\eta_{k+1} = 1\right\}} +  X^{NV,\eta}_{t_{k}} \mathds{1}_{\left\{\eta_{k+1} = -1\right\}}\right).
\end{equation}
Although the stochastic processes $\left(\bar{X}^{j,\eta}_{t}\right)_{t\in[0,T]}$, $j \in \left\{1,\ldots d\right\}$, are not adapted to the filtration $\mathbb{F}$, each stochastic integral is well defined in \eqref{NV-Interpol_ITO}. Indeed, $\left(\bar{X}^{j,\eta}_{t}\right)_{t\in[0,T]}$ is adapted with respect to the filtration $\left(\sigma\left(\eta, W^j_s, s\leq t   \right){\bigvee} \left(\underset{k \neq j}{\bigvee} \sigma\left( W^k_s, s\leq T   \right)\right)\right)_{t\in[0,T]}$, for $j \in \left\{1,\ldots d\right\}$. Then, by independence, $W^j$ is a also a Brownian motion with respect to this filtration and the stochastic integral $\displaystyle \int_0^t \displaystyle   \sigma^j(\bar{X}^{j,\eta}_s)  dW_s^j$ is well defined for all $t\in[0,T]$.
Using this interpolation, we proved in \cite{AJC} the strong convergence with order $1/2$. More precisely:
\begin{thm}
\label{SC_NV}
Assume that
\begin{itemize}
\item $\forall j \in \left\{1,\ldots,d\right\},\sigma^j \in \mathcal{C}^{1}\left(\mathbb{R}^n,\mathbb{R}^n\right)$.
\item $\sigma^0, \sigma^j$ and $\partial \sigma^j \sigma^j, \forall j \in \left\{1,\ldots,d\right\}$, are Lipschitz continuous functions.
\end{itemize}
Then, $\forall p \ge 1,  \exists C_{NV} \in \mathbb{R}^*_+,  \forall N \in \mathbb{N}^*$:
\begin{equation}
\label{SCG_NV}
\mathbb{E}\left[ \underset{t\leq T}{\sup}\left\|X_t - X^{NV,\eta}_{t}\right\|^{2p} \right] \leq C_{NV} \left(1+ \left\| x\right\|^{2p}\right) h^p.
\end{equation}
\end{thm} 
Then, the normalized error process is defined as follows
\begin{equation}
V^N = \sqrt{N}\left( X - X^{NV,\eta} \right).
\end{equation}
In this section, we check that the normalized error $V^N$ converges to an affine SDE with source terms. Here is the main result.
\begin{thm}
\label{EP_GC}
Assume that:
\begin{itemize}
\item  $\sigma^0 \in \mathcal{C}^2\left(\mathbb{R}^n,\mathbb{R}^n\right)$ and is a Lipschitz continuous function with polynomially growing second order derivatives.
\item  $\forall j \in \left\{1,\ldots,d\right\}, \sigma^j \in \mathcal{C}^2\left(\mathbb{R}^n,\mathbb{R}^n\right)$ and is Lipschitz continuous together with its first order derivative.
\item  $\forall j,m \in \left\{1,\ldots,d\right\},\partial \sigma^j\sigma^m$  is Lipschitz continuous.
\item  $\forall j \in \left\{1,\ldots,d\right\},\partial \sigma^j\sigma^j \in \mathcal{C}^2\left(\mathbb{R}^n,\mathbb{R}^n\right)$ with polynomially growing second order derivatives.

\end{itemize}
Then:
\begin{equation}
V^N \overset{stably}{\underset{N \to +\infty}{\Longrightarrow}} V
\end{equation}
where $V$ is the unique solution of the following affine equation:
\begin{equation}
V_t = \sqrt{\frac{T}{2}} \sum \limits_{j=1}^d \sum \limits_{m=1}^{j-1}  \int_0^t \displaystyle \left[\sigma^j,\sigma^m\right]\left(X_s\right) dB^{j,m}_s + \int_0^t \displaystyle \partial b\left(X_s \right)V_s ds + \sum \limits_{j=1}^d \int_0^t \displaystyle \partial \sigma^j\left(X_s \right)V_s dW_s^j
\label{EQ}
\end{equation} 
with $\left[\sigma^j,\sigma^m\right] = \partial \sigma^m \sigma^j - \partial \sigma^j \sigma^m $,
and $\left(B_t\right)_{0\leq t \leq T}$ a standard $\frac{d(d-1)}{2}$-dimensional Brownian motion independent of $W$.
\end{thm}

\subsection{Discrete scheme}

To compute the asymptotic error distribution, the method consists in writing the normalized error in the form \eqref{method_AED}
.
Since the interpolation \eqref{NV-Interpol_ITO} is not adapted to the natural filtration of the Brownian motion $W$, we were not able to derive a decomposition \eqref{method_AED} with $V^N$ replacing $U^N$. To get around this problem, we build an adapted approximation $\hat{X}^{D,\eta}$ of $X^{NV,\eta}$, with order $1-\epsilon, \forall \epsilon >0$, and introduce $U^N = \sqrt{N}\left(X - \hat{X}^{D,\eta}\right)$. Then, we obtain the decomposition of the form \eqref{method_AED} (see \eqref{Decomp_}) and study the satble convergence in law of $U^N$ to deduce the convergence of $V^N$.
The approximation is defined as follows:
\begin{equation}
\label{Approx_Discret}
\left\{
    \begin{array}{ll}
\begin{split}
\hat{X}^{D,\eta}_{t} & = \hat{X}^{D,\eta}_{\hat{\tau}_t} +  b\left(\hat{X}^{D,\eta}_{\hat{\tau}_t}\right) \Delta t +  \sum \limits_{j=1}^d \sigma^j\left(\hat{X}^{D,\eta}_{\hat{\tau}_t} \right) \Delta W^j_{t} + \frac{1}{2} \sum \limits_{j=1}^d \partial \sigma^j \sigma^j\left(\hat{X}^{D,\eta}_{\hat{\tau}_t}\right) \left( \left(\Delta W_{t}^j\right)^2 - \Delta t  \right)    \\
& +  \sum \limits_{\eta_t m < \eta_t j} \partial \sigma^j \sigma^m\left(\hat{X}^{D,\eta}_{\hat{\tau}_t}\right)\Delta W_{t}^m \Delta W_{t}^j     
\end{split}
\\
\hat{X}^{D,\eta}_0 = x.
\end{array}
\right.
\end{equation}

In the following proposition, we compare $X^{NV,\eta}$ and $\hat{X}^{D,\eta}$.
\begin{aprop}
\label{ADNV-NV}
Under the assumptions of Theorem \ref{EP_GC}:
\begin{equation}
\forall p \ge 1, \forall \epsilon > 0, \exists C_D \in \mathbb{R}_+^*, \forall N \in \mathbb{N}^*, ~ \mathbb{E}\left[ \underset{t\leq T}{\sup}\left\|X^{NV,\eta}_{t} - \hat{X}^{D,\eta}_{t}\right\|^{2p} \right] \leq C_D \frac{1}{N^{2p-\epsilon}}.
\end{equation}
\end{aprop}
The proof of this proposition is postponed to the Appendix.

The next lemma gives estimation of the moment of $\hat{X}^{D,\eta}$ and its increments. Its hypotheses are consequences of the ones of Theorem \ref{EP_GC}. We omit its standard proof.
\begin{alem}
\label{DNVA_Lemma}
Assume that:
\begin{itemize}
\item   $b\in \mathcal{C}^0\left(\mathbb{R}^n,\mathbb{R}^n\right)$ has an affine growth. 
\item $\forall j \in \left\{1,\ldots,d\right\}, \sigma^j$ has an affine growth. 
\item $\forall j,m \in \left\{1,\ldots,d\right\}, \partial \sigma^j\sigma^m$ has an affine growth. 
\end{itemize}
Then, $\forall p \ge 1, \exists \hat{C}_D \in \mathbb{R}_+^*, \forall N \in \mathbb{N}^*$:
 \begin{enumerate}[(i)]
 \item 
 \begin{equation} 	
 \label{DNVA_moment}
  \mathbb{E}\left[\underset{t\leq T}{\sup} \left\|\hat{X}^{D,\eta}_{t} \right\|^{2p} \right] \leq \hat{C}_D. 
 \end{equation}
\item 
\begin{equation}
\label{Increment}
\forall t \in [0,T], \mathbb{E}\left[ \left\|\hat{X}^{D,\eta}_{t} - \hat{X}^{D,\eta}_{\hat{\tau}_t}\right\|^{2p} \right] \leq \hat{C}_D h^p.
 \end{equation}
\end{enumerate}
\end{alem}

\subsection{Proof of the stable convergence}
We recall that $U^N = \sqrt{N}\left(X - \hat{X}^{D,\eta}\right)$.
By Proposition \ref{ADNV-NV}, $\underset{t\leq T}{\sup} \sqrt{N} \left\|\hat{X}^{D,\eta}_t - X^{NV,\eta}_t \right\|$ converges in probability to 0 as $N$ goes to $+\infty$. Since $V^N - U^N  = \sqrt{N}\left(\hat{X}^{D,\eta} - X^{NV,\eta} \right)$, \eqref{2ii} from Proposition \ref{Basic} ensures that Theorem \ref{EP_GC} is a consequence of the following proposition dedicated to the stable convergence in law of $U^N$.
\begin{aprop}
\label{Intermed}
Under the assumptions of Theorem \ref{EP_GC}:
\begin{equation}
U^N \overset{stably}{\underset{N \to +\infty}{\Longrightarrow}} V,
\end{equation}
where $V$ is the unique solution of \eqref{EQ}.
\end{aprop}
\begin{adem}
We begin by describing the limiting process for $U^N=\sqrt{N}\left( X - \hat{X}^{D,\eta} \right).$
The differential of $U^N$ can be written as:

\begin{equation}
\label{EP}
\begin{split}
dU^{N}_t & = \sqrt{N}  \left( \left(b\left(X_t\right) - b\left(\hat{X}^{D,\eta}_t\right) \right)dt  +  \sum \limits_{j=1}^d \left(\sigma^j\left(X_t\right) - \sigma^j\left(\hat{X}^{D,\eta}_t\right) \right)dW^j_t   \right)   \\
&+ \sqrt{N}  \left(\left(b\left(\hat{X}^{D,\eta}_t\right) - b\left(\hat{X}^{D,\eta}_{\hat{\tau}_t}\right) \right)dt  + \sum \limits_{j=1}^d \left(\sigma^j\left(\hat{X}^{D,\eta}_t\right) - \sigma^j\left(\hat{X}^{D,\eta}_{\hat{\tau}_t}\right) \right)dW^j_t   \right)  \\
& -\sqrt{N} \left(\sum \limits_{j=1}^d \partial \sigma^j \sigma^j\left(\hat{X}^{D,\eta}_{\hat{\tau}_t}\right) \Delta W_{t}^j dW_t^j ~ + \sum \limits_{\eta_t m < \eta_t j} \partial \sigma^j \sigma^m\left(\hat{X}^{D,\eta}_{\hat{\tau}_t}\right)\left( \Delta W_{t}^m dW_t^j  +\Delta W_{t}^j dW_t^m  \right) \right).
\end{split}
\end{equation}
Then, the proof will go through several steps. \\\\
\textbf{Step 1: linearisation of the two terms in the first line of the right-hand side of \eqref{EP}.}\\\\
\begin{equation}
 \sqrt{N} \left( \left(b\left(X_t\right) - b\left(\hat{X}^{D,\eta}_t\right) \right)dt  +  \sum \limits_{j=1}^d \left(\sigma^j\left(X_t\right) - \sigma^j\left(\hat{X}^{D,\eta}_t\right) \right)dW^j_t   \right) . 
\end{equation}
Let $j \in \left\{1,\ldots,d\right\}$ and $i \in \left\{1,\ldots,n\right\}$.
By the mean value theorem, we get:
\begin{equation}
 \sigma^{ij}\left(X_t\right) - \sigma^{ij}\left(\hat{X}^{D,\eta}_t \right) = \nabla \sigma^{ij}\left(\xi^{ij}_t\right) . \left(X_t - \hat{X}^{D,\eta}_t  \right)
\end{equation}
where:
$ \xi^{ij}_t = \alpha^{ij}_t X_t + \left(1 - \alpha^{ij}_t\right) \hat{X}^{D,\eta}_t $~ for some  $ \alpha^{ij}_t \in [0,1]$.
Using a compact matrix notation, we can write: 
\begin{equation}
 \sigma^{j}\left(X_t\right) - \sigma^{j}\left(\hat{X}^{D,\eta}_t\right) = \partial \sigma^{j,N}_t\left(X_t - \hat{X}^{D,\eta}_t \right)
\end{equation} 
where:
\begin{equation}
 \left(\partial \sigma^{j,N}_t\right)_{i,m} = \partial_{x_m} \sigma^{ij}\left(\xi_t^{ij}\right). 
\end{equation}
Then, we obtain
\begin{equation}
 \sqrt{N} \sum \limits_{j=1}^d \left(\sigma^j\left(X_t\right) - \sigma^j\left(\hat{X}^{D,\eta}_t\right) \right)dW^j_t  = \sum \limits_{j=1}^d \partial \sigma^{j,N}_t U^{N}_t dW^j_t.  
\end{equation}
In the same way:
\begin{equation}
 \sqrt{N} \left(b\left(X_t\right) - b\left(\hat{X}^{D,\eta}_t\right) \right)dt  = \partial b^N_t U^{N}_t dt  
\end{equation}
where:
\begin{equation}
 \left(\partial b^N_t\right)_{i,m} = \partial_{x_m} b^{i}\left(\xi_t^{i0}\right) 
\end{equation}
and $ \xi^{i0}_t = \alpha^{i0}_t X_t + \left(1 - \alpha^{i0}_t\right) \hat{X}^{D,\eta}_t $~ for some $ \alpha^{i0}_t \in [0,1]$.

\textbf{Step 2: decomposition of $U^N$.}\\\\
Writing the fourth term in the right-hand side of \eqref{EP}, $ \sigma^j\left(\hat{X}^{D,\eta}_t\right) - \sigma^j\left(\hat{X}^{D,\eta}_{\hat{\tau}_t}\right) $, as the sum of the dominant contribution $$\displaystyle \sum \limits_{m=1}^d \partial \sigma^j \sigma^{m}\left(\hat{X}^{D,\eta}_{\hat{\tau}_t}\right) \Delta W^m_t$$ and the remainder 
$$ \sigma^j\left(\hat{X}^{D,\eta}_t\right) - \sigma^j\left(\hat{X}^{D,\eta}_{\hat{\tau}_t}\right) - \displaystyle \sum \limits_{m=1}^d \partial \sigma^j \sigma^{m}\left(\hat{X}^{D,\eta}_{\hat{\tau}_t}\right) \Delta W^m_t,$$
which is of order $1/N$, we deduce that:
\begin{equation}
\label{Decomp_}
U^{N}_t = R^N_t + J^N_t + \left( \int_0^t \displaystyle \partial b_s^NU^N_s ds + \sum \limits_{j=1}^d \int_0^t \displaystyle \partial \sigma^{j,N}_s U^N_s dW_s^j \right)
\end{equation}
where
\begin{equation}
\begin{split}
 R^N_t &=   \sqrt{N} \Bigg( \sum \limits_{j=1}^d \int_0^t \displaystyle \left(  \sigma^j\left(\hat{X}^{D,\eta}_s\right) - \sigma^j\left(\hat{X}^{D,\eta}_{\hat{\tau}_s}\right) - \sum\limits_{m=1}^d \partial \sigma^j \sigma^m \left(\hat{X}^{D,\eta}_{\hat{\tau}_s}\right) \Delta W^m_s \right)  dW^j_s \\
&+ 	\int_0^t \displaystyle b\left(\hat{X}^{D,\eta}_s\right) - b\left(\hat{X}^{D,\eta}_{\hat{\tau}_s}\right) ds \Bigg)
\end{split}
 \end{equation} 
and
 \begin{equation}
\begin{split}
\label{JNt}
J^N_t & = -\sqrt{N} \Bigg(\sum \limits_{j=1}^d  \int_0^t \displaystyle  \partial \sigma^j \sigma^j\left(\hat{X}^{D,\eta}_{\hat{\tau}_s}\right) \Delta W_{s}^j dW_s^j ~ + \int_0^t \sum \limits_{\eta_s m < \eta_s j} \partial \sigma^j \sigma^m\left(\hat{X}^{D,\eta}_{\hat{\tau}_s}\right)\left\{ \Delta W_{s}^m dW_s^j  +\Delta W_{s}^j dW_s^m  \right\} \\
& - \sum \limits_{j=1}^d \sum \limits_{m=1}^d \int_0^t  \displaystyle \partial \sigma^j \sigma^{m}\left(\hat{X}^{D,\eta}_{\hat{\tau}_s}\right) \Delta W^m_s dW_s^j \Bigg).
\end{split}
\end{equation}
The expression of $J^N$ can be arranged as follows
\begin{equation}
\begin{split}
J^N_t & = \sqrt{N} \sum \limits_{j=1}^{d} \sum \limits_{m=1}^{j-1} \int_0^t \displaystyle \left[\sigma^j,\sigma^m\right] \left(\hat{X}^{D,\eta}_{\hat{\tau}_s}\right) \frac{-1 +\eta_s}{2}  \Delta W^m_s dW^j_s \\
& + \sqrt{N}  \sum \limits_{j=1}^{d} \sum \limits_{m=1}^{j-1} \int_0^t \displaystyle \left[\sigma^j,\sigma^m\right] \left(\hat{X}^{D,\eta}_{\hat{\tau}_s}\right) \frac{ 1 +\eta_s}{2}  \Delta W^j_s dW^m_s.
\end{split}
\end{equation}
\textbf{Step 3: stable convergence in law of $J^N$}\\\\
To lighten up the notations, we introduce:
\begin{itemize}
\item $K_t^{j,m,N} = \left[\sigma^j,\sigma^m\right]\left(\hat{X}^{D,\eta}_{\hat{\tau}_t}\right)$ for $j,m \left\{1,\ldots,d\right\}$, $m<j$.
\item $\Psi_t^1 =\frac{-1+\eta_t}{2} $ and $\Psi_t^2 = \frac{1+\eta_t}{2}$. 
\item $ Y^{j,m,N}_t = \sqrt{N} \left( \displaystyle  \int_0^t \displaystyle \Psi_s^1 \Delta W_s^m dW^j_s + \int_0^t \displaystyle \Psi_s^2 \Delta W_s^j dW^m_s \right)$ for $j,m \left\{1,\ldots,d\right\}$, $m<j$.
\end{itemize}
Then, we have that
\begin{equation}
\Psi^1 \Psi^2 = 0
\end{equation}
\begin{equation}
\left(\Psi^1\right)^2 + \left(\Psi^2\right)^2 = 1
\end{equation}
and
\begin{equation}
J^N_t = \sum \limits_{j=1}^{d} \sum \limits_{m=1}^{j-1} \displaystyle  \int_0^t \displaystyle K_s^{j,m,N} dY^{j,m,N}_s.
\end{equation}
With a view to apply Proposition \ref{Jacod_Protter}, in order to obtain the stable convergence in law of $J^N$, we first study the stable convergence in law of $Y^N$. By virtue of Theorem \ref{Jacod}, it is enough to study the asymptotic behavior of 
$ \left<Y^{j,m,N}, W^k \right>_t$, for $t \in [0,T], j,m,k \in \left\{1,\ldots,d\right\}$ such that $m < j$, and $\left<Y^{j,m,N}, Y^{l,k,N} \right>_t $, for $t \in [0,T], j,m,k,l \in \left\{1,\ldots,d\right\}$ such that $m < j$, $k <l$.

\textbf{Step 3.1: asymptotic behavior of $\left<Y^{j,m,N}, W^k \right>_t $ for $j,m,k \in \left\{1,\ldots,d\right\}, m < j$ and $t \in [0,T]$.}\\\\
Computing the quadratic covariation we get
\begin{equation}
\left<Y^{j,m,N}, W^k \right>_t = \sqrt{N} \left( \mathds{1}_{\left\{j=k\right\}} \int_0^t \displaystyle \Psi_s^1 \Delta W^m_s ds + \mathds{1}_{\left\{m=k\right\}} \int_0^t \displaystyle \Psi_s^2 \Delta W^j_s ds \right).
\end{equation}
Then, computing the $L^2-$norm, we obtain
\begin{equation}
\begin{split}
\left\| \left<Y^{j,m,N}, W^k \right>_t\right\|_2^2 & = N~ \mathbb{E} \left[\left( \mathds{1}_{\left\{j=k\right\}} \int_0^t \displaystyle \Psi_s^1 \Delta W^m_s ds + \mathds{1}_{\left\{m=k\right\}} \int_0^t \displaystyle \Psi_s^2 \Delta W^j_s ds \right)^2\right] \\
& = 2N \Bigg( \mathds{1}_{\left\{j=k\right\}}\int_0^t \displaystyle \int_0^s \displaystyle \mathbb{E} \left[\Psi_s^1 \Psi_u^1 \right] \mathbb{E} \left[ \Delta W^m_s \Delta W^m_u\right] du~ ds \\
&+ \mathds{1}_{\left\{m=k\right\}} \int_0^t \displaystyle \int_0^s \displaystyle \mathbb{E} \left[ \Psi_s^2 \Psi_u^2 \right]\mathbb{E} \left[ \Delta W^j_s \Delta W^j_u\right] du~ ds\Bigg). 
\end{split}
\end{equation}
Since for $u,s \in [0,t], u \leq s$, $\mathbb{E} \left[ \Delta W^m_s \Delta W^m_u\right] = u - u \wedge\hat{\tau}_s \ge 0$, $0 \leq \mathbb{E}\left[\Psi_s^1 \Psi_u^1 \right] = \mathbb{E}\left[\Psi_s^2 \Psi_u^2 \right] \leq \frac{1}{2}$ and $m<j$, then it follows that: 
\begin{equation}
\begin{split}
\left\| \left<Y^{j,m,N}, W^k \right>_t\right\|_2^2 & \leq N \int_0^t \displaystyle \int_0^s \displaystyle u - u \wedge\hat{\tau}_s ~ du~ ds = N \int_0^t \displaystyle \int_{\hat{\tau}_s}^s \displaystyle u - \hat{\tau}_s ~du~ ds = \frac12 N \int_0^t \displaystyle \left(s - \hat{\tau}_s\right)^2 ~ ds \\
& \leq  \frac12  N \int_0^T \displaystyle \left(s - \hat{\tau}_s\right)^2 ~ ds = \frac{T^3}{6N}  \underset{N \to +\infty}{\longrightarrow} 0.
\end{split}
\end{equation}
\textbf{Step 3.2: asymptotic behavior of $ \left<Y^{j,m,N}, Y^{l,k,N} \right>_t $ for $j,m,k,l \in \left\{1,\ldots,d\right\}, m < j, k <l $ and $t \in [0,T]$.}\\\\
Computing the quadratic covariation we get
\begin{equation}
\begin{split}
 \left<Y^{j,m,N}, Y^{l,k,N} \right>_t & = N  \Bigg( \mathds{1}_{\left\{j=l\right\}}   \int_0^t \displaystyle  \left(\Psi_s^1\right)^2 \Delta W^m_s \Delta W^k_s ds  + \mathds{1}_{\left\{m=k\right\}}   \int_0^t \displaystyle  \left(\Psi_s^2\right)^2 \Delta W^j_s \Delta W^l_s ds  \\
& + \mathds{1}_{\left\{j=k\right\}}   \int_0^t \displaystyle \Psi_s^1 \Psi_s^2 \Delta W^l_s \Delta W^m_s ds +  \mathds{1}_{\left\{m=l\right\}}   \int_0^t \displaystyle \Psi_s^1 \Psi_s^2 \Delta W^j_s \Delta W^k_s ds  \Bigg).
\end{split}
\end{equation}
Since $\Psi^1 \Psi^2 = 0$, we obtain 
\begin{equation}
\begin{split}
 \left<Y^{j,m,N}, Y^{l,k,N} \right>_t = N  \left( \mathds{1}_{\left\{j=l\right\}} \int_0^t \displaystyle  \left(\Psi_s^1\right)^2 \Delta W^m_s \Delta W^k_s ds  + \mathds{1}_{\left\{m=k\right\}}  \int_0^t \displaystyle  \left(\Psi_s^2\right)^2 \Delta W^j_s \Delta W^l_s ds \right).
\end{split}
\end{equation}
Then, we distinguish three cases. In the case $j = l$ and $m \neq k$, the bracket is given by
\begin{equation}
\begin{split}
 \left<Y^{j,m,N}, Y^{l,k,N} \right>_t = N   \int_0^t \displaystyle  \left(\Psi_s^1\right)^2 \Delta W^m_s \Delta W^k_s ds.
\end{split}
\end{equation}
By independence, the $L^2-$norm of the bracket $ \left<Y^{j,m}, Y^{j,k} \right>_t $ is given by 
\begin{equation}
\begin{split}
\left\|  \left<Y^{j,m,N}, Y^{j,k,N} \right>_t \right\|_2^2 & = 2N^2 \int_0^t \displaystyle \int_0^s \displaystyle \mathbb{E} \left[\left(\Psi_s^1\right)^2 \left(\Psi_u^1\right)^2 \right]\mathbb{E} \left[ \Delta W^m_s \Delta W^m_u\right] \mathbb{E} \left[ \Delta W^k_s \Delta W^k_u\right] du~ ds \\
 & =  2N^2 \int_0^t \displaystyle \int_0^s \displaystyle \mathbb{E} \left[\left(\Psi_s^1\right)^2 \left(\Psi_u^1\right)^2 \right]\left(\mathbb{E} \left[ \Delta W^m_s \Delta W^m_u\right]\right)^2 du~ ds. 
\end{split}
\end{equation}
Once again, using $\mathbb{E} \left[ \Delta W^m_s \Delta W^m_u\right] = u - u \wedge\hat{\tau}_s \ge 0$, and $0\leq \mathbb{E}\left[\left(\Psi_s^1\right)^2 \left(\Psi_u^1\right)^2 \right]  \leq \frac{1}{2}$ for $u,s \in [0,t], u \leq s$, we get
\begin{equation}
\begin{split}
\left\|  \left<Y^{j,m,N}, Y^{j,k,N} \right>_t \right\|_2^2 & \leq N^2 \int_0^t \displaystyle \int_0^s \displaystyle \left(u - u \wedge\hat{\tau}_s\right)^2 du~ ds  =  N^2 \int_0^t \displaystyle \int_{\hat{\tau}_s}^s \displaystyle \left(u -\hat{\tau}_s\right)^2 du~ ds    = \frac 13 N^2 \int_0^t \displaystyle \displaystyle \left(s -\hat{\tau}_s\right)^3  ds\\
& \leq \frac 13 N^2 \int_0^T \displaystyle \displaystyle \left(s -\hat{\tau}_s\right)^3  ds = \frac{1}{12} \frac{T^4}{N}  \underset{N \to +\infty}{\longrightarrow} 0.
\end{split}
\end{equation}
The second case $k = m$ and $j \neq l$, is similar to the first case since
\begin{equation}
 \left<Y^{j,k,N}, Y^{l,k,N} \right>_t = N  \int_0^t \displaystyle \left(\Psi_s^2\right)^2 \Delta W^j_s \Delta W^l_s ds.   
\end{equation}
As previously, we have:
\begin{equation}
\left\|  \left<Y^{j,k,N}, Y^{l,k,N} \right>_t  \right\|_2 \underset{N \to +\infty}{\longrightarrow} 0.
\end{equation}
The third and last case $k = m$ and $j = l$ provides a nonzero limit
\begin{equation}
 \left<Y^{j,m,N}, Y^{j,m,N} \right>_t = N  \int_0^t \displaystyle \left(\left(\Psi_s^1\right)^2\left( \Delta W^m_s \right)^2 + \left(\Psi_s^2\right)^2\left( \Delta W^j_s \right)^2 \right) ds.   
\end{equation}
To identify the limit we proceed to a preliminary calculus of the expectation of this bracket
\begin{equation}
\mathbb{E}\left[ \left<Y^{j,m,N}, Y^{j,m,N} \right>_t \right]= N  \int_0^t \displaystyle \left(\mathbb{E} \left[\left(\Psi_s^1\right)^2\right] \mathbb{E}\left[ \left( \Delta W^m_s \right)^2\right] + \mathbb{E} \left[\left(\Psi_s^2\right)^2\right]\mathbb{E}\left[ \left( \Delta W^j_s \right)^2 \right] \right)ds.   
\end{equation}
Since $\left(\Psi^1\right)^2 +  \left(\Psi^2\right)^2 = 1$, we get
\begin{equation}
\begin{split}
\mathbb{E}\left[ \left<Y^{j,m,N}, Y^{j,m,N} \right>_t \right]& = N  \int_0^t \displaystyle \left(s - \hat{\tau}_s\right)~ ds   = N \int_0^{\hat{\tau}_t} \displaystyle \left(s - \hat{\tau}_s\right)~ ds   + O\left(\frac{1}{N}\right) = \frac{1}{2} Tt + O\left(\frac{1}{N}\right) \underset{N \to +\infty}{\longrightarrow} \frac{1}{2} Tt .
\end{split}
\end{equation}
Now, we show the convergence in $L^2$  of $\left<Y^{j,m,N}, Y^{j,m,N} \right>_t$  towards $\frac{1}{2} Tt$. Computing the $L^2-$norm of the difference between the quadratic variation $\left<Y^{j,m}, Y^{j,m} \right>_t$ and $\frac{1}{2} Tt $, we obtain
\begin{equation}
\begin{split}
\left\|  \left<Y^{j,m,N}, Y^{j,m,N} \right>_t  - \frac{1}{2} Tt \right\|^2_2 & = \left\|  \left<Y^{j,m,N}, Y^{j,m,N} \right>_t \right\|^2_2  -2 \mathbb{E}\left[ \left<Y^{j,m,N}, Y^{j,m,N} \right>_t \right] \frac{1}{2} Tt + \left(\frac{1}{2} Tt\right)^2\\
&  = \left\|  \left<Y^{j,m,N}, Y^{j,m,N} \right>_t \right\|^2_2  -  \left(\frac{1}{2} Tt\right)^2  + O\left(\frac{1}{N}\right).
\end{split}
\end{equation}
To prove the convergence in $L^2$, it suffices to show that
\begin{equation}
 \left\|  \left<Y^{j,m,N}, Y^{j,m,N} \right>_t \right\|^2_2  \underset{N \to +\infty}{\longrightarrow} \left(\frac{1}{2} Tt \right)^2.
\end{equation}
 Computing the square of the $L^2$ norm of the bracket
\begin{equation}
\begin{split}
\left\|  \left<Y^{j,m,N}, Y^{j,m,N} \right>_t  \right\|_2^2 &= 2N^2\Bigg(  \int_0^t \displaystyle \int_0^s \displaystyle \mathbb{E} \left[\left(\Psi_s^1\right)^2 \left(\Psi_u^1\right)^2\right] \mathbb{E}\left[ \left( \Delta W^m_s \right)^2\left( \Delta W^m_u \right)^2\right] du~ ds\\
 &+\int_0^t \displaystyle \int_0^s \displaystyle \mathbb{E} \left[ \left(\Psi_s^1\right)^2 \left(\Psi_u^2\right)^2\right] \mathbb{E}\left[ \left( \Delta W^m_s \right)^2\left( \Delta W^j_u \right)^2\right] du~ ds\\
 &+\int_0^t \displaystyle \int_0^s \displaystyle \mathbb{E} \left[ \left(\Psi_s^2\right)^2 \left(\Psi_u^1\right)^2\right] \mathbb{E}\left[ \left( \Delta W^j_s \right)^2\left( \Delta W^m_u \right)^2\right] du~ ds\\
 &+\int_0^t \displaystyle \int_0^s \displaystyle \mathbb{E} \left[\left(\Psi_s^2\right)^2 \left(\Psi_u^2\right)^2\right]\mathbb{E}\left[ \left( \Delta W^j_s \right)^2\left( \Delta W^j_u \right)^2\right] du~ ds\Bigg).
\end{split} 
\end{equation}
Since for $k,l \in \left\{1,\ldots,d\right\}$, $\mathbb{E}\left[ \left( \Delta W^k_s \right)^2\left( \Delta W^l_u \right)^2\right] = O\left(\frac{1}{N^2}\right) , \forall u,s \in [0,t], u \leq s$, we get
\begin{equation}
\begin{split}
\left\|  \left<Y^{j,m,N}, Y^{j,m,N} \right>_t  \right\|_2^2 &= 2N^2 \int_0^t \displaystyle \int_0^{\hat{\tau}_s}  \displaystyle \mathbb{E} \left[ \left( \left(\Psi_u^1\right)^2 +  \left(\Psi_u^2\right)^2 \right)\left( \left(\Psi_s^1\right)^2 +  \left(\Psi_s^2\right)^2 \right)  \right] \left(u - \hat{\tau}_u\right)  \left(s - \hat{\tau}_s\right)  du ~ds +  O\left(\frac{1}{N}\right).
\end{split} 
\end{equation}
Then, using $\left(\Psi^1\right)^2 + \left(\Psi^2\right)^2 = 1$
\begin{equation}
\begin{split}
\left\|  \left<Y^{j,m,N}, Y^{j,m,N} \right>_t  \right\|_2^2 
& = 2 N^2 \int_0^t \displaystyle \int_0^{\hat{\tau}_s}  \displaystyle  \left(u - \hat{\tau}_u\right)  \left(s - \hat{\tau}_s\right)  du ~ds +  O\left(\frac{1}{N}\right)\\
& = 2 N^2 \int_0^{t} \displaystyle \int_0^{s}  \displaystyle  \left(u - \hat{\tau}_u\right)  \left(s - \hat{\tau}_s\right)  du ~ds +  O\left(\frac{1}{N}\right)\\
& =  N^2 \left(\int_0^{t}  \left(s - \hat{\tau}_s\right) ds \right)^2 +  O\left(\frac{1}{N}\right) = \left(\frac{1}{2} Tt \right)^2 +  O\left(\frac{1}{N}\right)  \underset{N \to +\infty}{\longrightarrow} \left(\frac{1}{2} Tt \right)^2.
\end{split} 
\end{equation}
\textbf{Step 3.3: conclusion of the step 3.}\\\\
Applying  Theorem \ref{Jacod} we conclude that $\sqrt{\frac{2}{T}}Y^N$ stably converges in law to a standard $\frac{d(d-1)}{2}$-dimensional Brownian motion $B$, independent of $W$.

Now, it remains to prove the convergence in probability of $K^N$. We recall that for $j,m \left\{1,\ldots,d\right\}$, $m<j$, $K_t^{j,m,N} = \left[\sigma^j,\sigma^m\right]\left(\hat{X}^{D,\eta}_{\hat{\tau}_t}\right)$ . From Proposition \ref{ADNV-NV} and Theorem \ref{SC_NV}, together with the Lipschitz assumption on $\left[\sigma^j,\sigma^m\right]$, $j,m \left\{1,\ldots,d\right\}$, $m<j$ we get the following convergence in $L^2$
\begin{equation}
\label{CV_KL2}
\underset{t\leq T}{\sup} \left\| K_t^{j,m,N} - \left[\sigma^j,\sigma^m\right]\left(X_{\hat{\tau}_t}\right) \right\| \overset{L^2}{\underset{N \to +\infty}{\longrightarrow}} 0.
\end{equation}
Once again, the continuity of $\left[\sigma^j,\sigma^m\right]$, $j,m \left\{1,\ldots,d\right\}$, $m<j$, together with the continuity of the solution $X$ ensure that 
\begin{equation}
\label{CV_KP}
\underset{t\leq T}{\sup} \left\| \left[\sigma^j,\sigma^m\right]\left(X_{\hat{\tau}_t}\right) - \left[\sigma^j,\sigma^m\right]\left(X_{t}\right) \right\| \overset{a.s}{\underset{N \to +\infty}{\longrightarrow}} 0.
\end{equation}
Then, combining \eqref{CV_KL2} and \eqref{CV_KP}, we obtain
\begin{equation}
\underset{t\leq T}{\sup} \left\| K_t^{j,m,N} - \left[\sigma^j,\sigma^m\right]\left(X_{t}\right) \right\| \overset{\mathbb{P}}{\underset{N \to +\infty}{\longrightarrow}} 0.
\end{equation}
Thus, according to Proposition \ref{Basic}, we have the following convergence:
\begin{equation}
\left(K^N, \sqrt{\frac{2}{T}}Y^N\right) \overset{stably}{\underset{N \to +\infty}{\Longrightarrow}} \left(\left(\left[\sigma^j,\sigma^m\right]\left(X\right)\right)_{j,m \in \left\{1,\ldots,d\right\}, m<j},B\right)
\end{equation}
The convergence of $\left<Y^N\right>_T$ ensures its tightness. Then Proposition \ref{Jacod_Protter} leads us to:
\begin{equation}
\label{Step3}
\left(K^N, \sqrt{\frac{2}{T}}Y^N, J^N\right) \overset{stably}{\underset{N \to +\infty}{\Longrightarrow}} \left(\left(\left[\sigma^j,\sigma^m\right]\left(X\right)\right)_{j,m \in \left\{1,\ldots,d\right\}, m<j},B,\left(\sqrt{\frac{T}{2}} \sum \limits_{j=1}^d \sum \limits_{m=1}^{j-1}  \int_0^t \displaystyle \left[\sigma^j,\sigma^m\right]\left(X_s\right) dB^{j,m}_s \right)_{t\in [0,T]} \right).
\end{equation}
\textbf{Step 4: convergence of $R^N$.}\\\\
We show the convergence in $L^2$ of the remainder $R^N$ towards 0.
Applying a convexity inequality, then Doob's martingale inequality to each stochastic integral, we get
\begin{equation}
\mathbb{E} \left[\underset{t\leq T}{\sup}\left\|R^N_t\right\|^{2}\right] \leq N \left(d+1\right) \left( E_0 + 4 \sum \limits_{j=1}^d E_j\right)
\end{equation}
where:
\begin{equation}
E_0 =  \mathbb{E} \left[\underset{t\leq T}{\sup}\left\| \int_0^t \displaystyle b\left(\hat{X}^{D,\eta}_s\right) - b\left(\hat{X}^{D,\eta}_{\hat{\tau}_s}\right) ds \right\|^{2}\right]
\end{equation}
and, for $j\in \left\{1,\ldots,d\right\}$
\begin{equation}
E_j = \int_0^T \displaystyle  \mathbb{E} \left[ \left\| \sigma^j\left(\hat{X}^{D,\eta}_s\right) - \sigma^j\left(\hat{X}^{D,\eta}_{\hat{\tau}_s}\right) - \sum\limits_{m=1}^d \partial \sigma^j \sigma^m \left(\hat{X}^{D,\eta}_{\hat{\tau}_s}\right) \Delta W^m_s \right\|^{2}   \right] ds. 
\end{equation}
\textbf{Step 4.1: estimation of $E_0$.}\\\\
Let $i\in \left\{1,\ldots,n\right\}$. Using the integration by parts and It\^o's formulae we get:
\begin{equation}
\begin{split}
\int_0^t \displaystyle b^i\left(\hat{X}^{D,\eta}_s\right) - b^i\left(\hat{X}^{D,\eta}_{\hat{\tau}_s}\right) du &= \int_0^t \displaystyle \left( t \wedge \check{\tau}_s - s \right) \nabla b^i\left(\hat{X}^{D,\eta}_s\right) .~ b^i\left(\hat{X}^{D,\eta}_{\hat{\tau}_s}\right) ds \\
&+ \sum \limits_{j=1}^d \int_0^t \displaystyle \left( t \wedge \check{\tau}_s - s \right) \nabla b^i\left(\hat{X}^{D,\eta}_s\right) .~ \sigma^j\left(\hat{X}^{D,\eta}_{\hat{\tau}_s}\right) dW^j_s\\
&+ \sum \limits_{j=1}^d \int_0^t \displaystyle \left( t \wedge \check{\tau}_s - s \right) \nabla b^i\left(\hat{X}^{D,\eta}_s\right) .~ \partial \sigma^j \sigma^j\left(\hat{X}^{D,\eta}_{\hat{\tau}_s}\right) \Delta W^j_s dW^j_s\\
&+ \int_0^t \displaystyle \sum \limits_{\eta_s m < \eta_s j }\left( t \wedge \check{\tau}_s - s \right) \nabla b^i\left(\hat{X}^{D,\eta}_s\right) .~ \partial \sigma^j \sigma^m\left(\hat{X}^{D,\eta}_{\hat{\tau}_s}\right) \Delta W^m_s dW^j_s\\
&+ \int_0^t \displaystyle \sum \limits_{\eta_s m < \eta_s j }\left( t \wedge \check{\tau}_s - s \right) \nabla b^i\left(\hat{X}^{D,\eta}_s\right) .~ \partial \sigma^j \sigma^m\left(\hat{X}^{D,\eta}_{\hat{\tau}_s}\right) \Delta W^j_s dW^m_s\\
&+ \frac12 \int_0^t \displaystyle \left( t \wedge \check{\tau}_s - s \right) Tr\left(H_s^* H_s \nabla^2 b^i\left(\hat{X}^{D,\eta}_s\right)\right )ds\\
\end{split}
\end{equation}
where $\left(H_s\right)_{0\leq s \leq T}$ is a $n\times d$-dimensional process built with the following columns:
\begin{equation*}
H^j_s = \sigma^j\left(\hat{X}^{D,\eta}_{\hat{\tau}_s}\right) + \partial \sigma^j \sigma^j  \left(\hat{X}^{D,\eta}_{\hat{\tau}_s}\right)\Delta W^j_s  + \sum \limits_{\eta_s m < \eta_s j } \partial \sigma^j \sigma^m  \left(\hat{X}^{D,\eta}_{\hat{\tau}_s}\right) \Delta W^m_s + \sum \limits_{\eta_s m > \eta_s j } \partial \sigma^m \sigma^j  \left(\hat{X}^{D,\eta}_{\hat{\tau}_s}\right) \Delta W^m_s. 
\end{equation*}
Since the first and second order derivatives of $b = \sigma^0 + \frac{1}{2} \sum \limits_{j=1}^d \partial \sigma^j \sigma^j$ have a polynomial growth and $t \wedge \check{\tau}_s - s \leq \frac{T}{N} $, \eqref{DNVA_moment} from Lemma \eqref{DNVA_Lemma} ensures the existence of a constant $\beta \in \mathbb{R}_+^*$ independent of $N$ such that:
\begin{equation}
\mathbb{E} \left[\underset{t\leq T}{\sup}\left| \int_0^t \displaystyle b^i\left(\hat{X}^{D,\eta}_s\right) - b^i\left(\hat{X}^{D,\eta}_{\hat{\tau}_s}\right) ds \right|^{2}\right] \leq\frac{\beta}{N^2}.
\end{equation}
\textbf{Step 4.2: estimation of $E_j, j \in \left\{1,\ldots,d\right\}$.}\\\\
Let $i\in \left\{1,\ldots,n\right\}$. Denoting
\begin{equation}
\begin{split}
\Phi^{ij}_s =\sigma^{ij}\left(\hat{X}^{D,\eta}_s\right) - \sigma^{ij}\left(\hat{X}^{D,\eta}_{\hat{\tau}_s}\right) - \sum\limits_{m=1}^d \nabla \sigma^{ij}\left(\hat{X}^{D,\eta}_{\hat{\tau}_s}\right).~ \sigma^m \left(\hat{X}^{D,\eta}_{\hat{\tau}_s}\right) \Delta W^m_s 
\end{split}
\end{equation}
and applying the mean value theorem, we get a $\zeta^{ij}_s$ between $\hat{X}^{D,\eta}_s$ and $\hat{X}^{D,\eta}_{\hat{\tau}_s}$, such that:
\begin{equation}
\begin{split}
\Phi^{ij}_s &= \sum \limits_{m=1}^d   \left(\nabla \sigma^{ij}\left(\zeta^{ij}_s\right) -  \nabla \sigma^{ij}\left(\hat{X}^{D,\eta}_{\hat{\tau}_s}\right) \right).~ \sigma^m\left(\hat{X}^{D,\eta}_{\hat{\tau}_s}\right)  \Delta W^m_s + \nabla \sigma^{ij}\left(\zeta^{ij}_s\right).~\displaystyle\int_{\hat{\tau}_s}^{s} b\left(\hat{X}^{D,\eta}_{\hat{\tau}_u}\right) du  \\
&+  \nabla \sigma^{ij}\left(\zeta^{ij}_s\right).~ \Bigg( \sum \limits_{\eta_s k < \eta_s m} \displaystyle\int_{\hat{\tau}_s}^{s}  \partial \sigma^m \sigma^k\left(\hat{X}^{D,\eta}_{\hat{\tau}_u}\right)\left\{ \Delta W_{u}^k dW_u^m  +\Delta W_{u}^m dW_u^k  \right\}\\
& + \sum \limits_{m=1}^d  \displaystyle\int_{\hat{\tau}_s}^{s} \partial \sigma^m \sigma^m\left(\hat{X}^{D,\eta}_{\hat{\tau}_u}\right) \Delta W_{u}^m dW_u^m  \Bigg).
\end{split}
\end{equation}
Since $b,\sigma^j, \partial \sigma^j,\partial \sigma^j \sigma^m, j,m \in \left\{1,\ldots,d\right\},$ are Lipschitz continuous, using  \eqref{DNVA_moment} and \eqref{Increment} from Lemma \ref{DNVA_Lemma}, we get a constant $\gamma\in \mathbb{R}_+$ independent of $N$ such that:
\begin{equation}
\mathbb{E} \left[ \left| \Phi^{ij}_s \right|^2\right] \leq\frac{\gamma}{N^2}.
\end{equation}
Then, it follows that:
\begin{equation}
\int_0^T \displaystyle  \mathbb{E} \left[ \left| \sigma^{ij}\left(\hat{X}^{D,\eta}_s\right) - \sigma^{ij}\left(\hat{X}^{D,\eta}_{\hat{\tau}_s}\right) - \sum\limits_{m=1}^d \nabla \sigma^{ij}\left(\hat{X}^{D,\eta}_{\hat{\tau}_s}\right).~ \sigma^m \left(\hat{X}^{D,\eta}_{\hat{\tau}_s}\right) \Delta W^m_s \right|^{2}   \right] ds \leq \frac{\gamma T}{N^2}.
\end{equation}
\textbf{Step 4.3: conclusion of the step 4.}\\\\
Taking the Euclidean norm, we conclude that:
\begin{equation}
\label{Step4}
\mathbb{E} \left[\underset{t\leq T}{\sup}\left\|R^N_t\right\|^{2}\right] \underset{N \rightarrow +\infty }{\longrightarrow} 0. 
\end{equation}
\textbf{Step 5: stable convergence in law of $U^N$.}\\\\
We recall that $U^{N}_t = R^N_t + J^N_t + \left( \displaystyle \int_0^t  \partial b_s^NU^N_s ds + \sum \limits_{j=1}^d \displaystyle \int_0^t \displaystyle \partial \sigma^{j,N}_s U^N_s dW_s^j \right)$.
Thanks to \eqref{Step3} and \eqref{Step4}, we conclude using Theorem \ref{Jacod_Protter_EDS} since the continuity of $\partial b$ and $\partial \sigma^j, j \in  \left\{1,\ldots,d\right\},$ together with Proposition \ref{ADNV-NV} ands Theorem \ref{SC_NV}, ensure that 
\begin{equation}
\underset{t\leq T}{\sup} \left\| \partial b^N_t - \partial b\left(X_t\right) \right\| \overset{\mathbb{P}}{\underset{N \to +\infty}{\longrightarrow}} 0,
\end{equation}
and $\forall j \in \left\{1,\ldots,d\right\}$,
\begin{equation}
\underset{s\leq T}{\sup} \left\| \partial \sigma^{j,N}_t - \partial \sigma^{j}\left(X_t\right) \right\| \overset{\mathbb{P}}{\underset{N \to +\infty}{\longrightarrow}} 0.
\end{equation}
\end{adem}

\section{Particular case: the Brownian vector fields commute}
In this section, we assume the following commutativity condition $$\forall j,m \in \left\{1,\ldots,d\right\}, \left[ \sigma^j,\sigma^m\right] = \partial \sigma^m\sigma^j - \partial \sigma^j \sigma^m = 0. $$ The commutativity of the Brownian vector fields implies the commutativity of the associated flows. Then, the order of integration of these fields no longer matters and $\eta$ is useless.
We also assume the following regularity assumptions:
\begin{itemize}
\item $\forall j \in \left\{1,\ldots,d\right\}, \sigma^j \in \mathcal{C}^{1}\left(\mathbb{R}^n,\mathbb{R}^n\right)$  with bounded first order derivatives.
\item $\sigma^0 \in \mathcal{C}^{2}\left(\mathbb{R}^n,\mathbb{R}^n\right)$ with bounded first order derivatives and polynomially growing second order derivatives.
\item $\sum \limits_{j=1}^d \partial \sigma^j \sigma^j$ is a Lipschitz continuous function.
\end{itemize}
Notice that $b = \sigma^0 - \frac12 \sum \limits_{j=1}^d \partial \sigma^j \sigma^j$ is also Lipschitz continuous.
We denote by $L \in \mathbb{R}_+^*$ a common Lipschitz constant of $\sigma^j, j \in \left\{0,\ldots,d\right\}$, $b$ and  $\sum \limits_{j=1}^d \partial \sigma^j \sigma^j$.
When the vector fields corresponding to each Brownian coordinate in the SDE commute, the solution of \eqref{EQ} is zero. This suggests that the rate of convergence  
is greater than $1/2$ in this case. In fact, we prove strong convergence with order $1$.

\subsection{Interpolated scheme and strong convergence}
In the commutative case we can define a smart interpolation. We define three intermediate processes. For $t \in \left(t_k, t_{k+1}\right]$:
\begin{equation}
\bar{X}^0_t = \exp\left(\frac{\Delta t}{2}\sigma^0\right)  X^{NV}_{t_{k}},
\end{equation}
\begin{equation}
\bar{X}_t = \exp\left(\Delta W_t^d\sigma^d\right) \ldots \exp\left(\Delta W_t^1\sigma^1\right)   \bar{X}^{0}_{t_{k+1}},
\label{XBar}
\end{equation}
\begin{equation}
\bar{X}^{d+1}_t = \exp\left(\frac{\Delta t}{2}\sigma^0\right) \bar{X}_{t_{k+1}}, 
\end{equation}
\begin{equation}
 X^{NV}_{t_{k+1}} = \bar{X}^{d+1}_{t_{k+1}}. 
\end{equation}
\begin{aprop}
\label{HD}
Let $t \in \left(t_k, t_{k+1}\right]$. The dynamics of $\left( \bar{X}_t\right)_{t \in \left(t_k, t_{k+1}\right]}$ is given by:
\begin{equation}
\bar{X}_t = \bar{X}^{0}_{t_{k+1}} +  \sum \limits_{j=1}^d \displaystyle\int_{t_k}^{t} \sigma^j \left(\bar{X}_s\right) \circ dW^j_s.
\end{equation} 
\end{aprop}
\begin{adem}
Frobenius' theorem (see \cite{Doss} or \cite{Dieudonee}) ensures that there exists a unique function $h \in \mathcal{C}^{1,2}\left(\mathbb{R}^n\times\mathbb{R}^d,\mathbb{R}^n\right)$ such that: 
\begin{equation}
\left\{
    \begin{array}{ll}
\partial_y h(x,y) = \sigma\left(( h(x,y)\right) ~~ \forall (x,y) \in \mathbb{R}^n \times \mathbb{R}^d\\\
h(x,0) = x ~~ \forall x \in \mathbb{R}^n
\end{array}
\right.
\label{Frob}
\end{equation}
ie $\forall j \in \left\{1,\ldots,d\right\}$
\begin{equation*}
\partial_{y_j} h(x,y) = \sigma^j\left(( h(x,y)\right) ~~ \forall (x,y) \in \mathbb{R}^n \times \mathbb{R}^d.
\end{equation*}
Then, it is clear that, by induction on $j \in \left\{1,\ldots,d\right\}$:
\begin{equation*}
 \exp\left(\Delta W_t^j\sigma^j\right) \ldots \exp\left(\Delta W_t^1\sigma^1\right)   \bar{X}^{0}_{t_{k+1}} =h\left(\bar{X}^{0}_{t_{k+1}} ,\Delta W^1_t,\ldots,\Delta W^j_t,0,\ldots,0\right) ,
\end{equation*}
and
\begin{equation*}
\bar{X}_t = h\left(\bar{X}^{0}_{t_{k+1}} ,\Delta W_t\right). 
\end{equation*}
 Finally, we obtain the desired result by applying the chain rule for the Stratonovich integral.
\end{adem}
The interpolated scheme is defined as follows:
\begin{equation}
X^{NV}_{t} = x + \frac{1}{2} \displaystyle\int_{0}^{t} \sigma^0 \left(\bar{X}^0_s\right)  ds +  \sum \limits_{j=1}^d \displaystyle\int_{0}^{t} \sigma^j \left(\bar{X}_s\right) \circ dW^j_s +  \frac{1}{2} \displaystyle\int_{0}^{t} \sigma^0 \left(\bar{X}^{d+1}_s\right) ds.
\end{equation} 

\begin{thm}
\label{SC2}
Let $p \in [1,+\infty)$.
Under the commutativity assumption, and the regularity assumptions made in the beginning of the section, we have the following convergence rate
\begin{equation}
\exists C^{\prime}_{NV}\in \mathbb{R}_+^*, \forall N \in \mathbb{N}^*, ~ \mathbb{E}\left[ \underset{t\leq T}{\sup}\left\|X_t - X^{NV}_{t}\right\|^{2p} \right] \leq C^{\prime}_{NV} h^{2p}.
\end{equation}
\end{thm}
When all the vector fields commute, the Ninomiya-Victoir scheme solves the SDE \eqref{EDS_ITO}. This suggests that the asymptotic distribution of the normalized error process $N\left(X-X^{NV}\right)$ involves the Lie brackets between the Brownian vector fields and the drift vector field.
We plan to investigate this question in a further work.
In order to prove Theorem \ref{SC2}, we first need to prove that the Ninomiya-Victoir scheme has uniformly bounded moments under the assumptions made in the beginning of this section. This is the aim of the following subsection.
\subsection{Intermediate results}

\begin{aprop}
\label{P1}
Let $p \ge 1$, $Z= \left( Z_t \right)_{0\leq t \leq h}$ and $Y= \left( Y_t \right)_{0\leq t \leq h}$ be the solutions of the following $n-$dimensional SDEs, driven by a $d-$dimensional brownian motion, on the time interval $[0,h]$:
\begin{equation*}
\left\{
    \begin{array}{ll}
 dZ_s  = \alpha(Z_s) ds + \beta(Z_s) dW_s \\
  Z_0 \text{ independent of } \left(W_t\right)_{t \in [0,h]} \text{such that:~} \mathbb{E}\left[\left\| Z_0 \right\|^{2p} \right] <+ \infty,
\end{array}
\right.
\end{equation*}
\begin{equation*}
\left\{
    \begin{array}{ll}
 dY_s  = \gamma(Y_s) ds + \beta(Y_s) dW_s \\
  Y_0  \text{ independent of } \left(W_t\right)_{t \in [0,h]} \text{such that:~} \mathbb{E}\left[\left\| Y_0 \right\|^{2p} \right] <+ \infty,
\end{array}
\right.
\end{equation*}
respectively. Assume that $\alpha$, $\beta$ and $\gamma$ are Lipschitz continuous functions. Then, $\exists C_0 \in \mathbb{R}^*_+, \forall t,s \in [0,h], s \leq t,$
\begin{enumerate}[(i)]
\item \begin{equation}
 \mathbb{E}\left[ 1 + \left\| Z_t \right\|^{2p} \right] \leq \mathbb{E}\left[ 1 + \left\| Z_0 \right\|^{2p} \right] \exp\left(C_0h\right). 
\label{R1}
\end{equation}
\item \begin{equation}
 \mathbb{E}\left[ \underset{t\leq h}{\sup}\left\| Z_t - Y_t \right\|^{2p}  \right] \leq C_0 \left( \mathbb{E}\left[ \left\| Z_0 - Y_0 \right\|^{2p} \right] + \left( \mathbb{E}\left[1 + \left\| Y_0 \right\|^{2p}\right] + \mathbb{E}\left[1 + \left\| Z_0 \right\|^{2p}\right] \right) h^{2p}\right) \exp\left(C_0 h \right).
\label{R2}
\end{equation}
If $\alpha = \gamma$, we have:
\begin{equation}
\label{R3}
\mathbb{E}\left[ \left\| Z_t - Y_t \right\|^{2p} \right] \leq \mathbb{E}\left[ \left\| Z_0 - Y_0 \right\|^{2p}\right] \exp \left(C_0 h \right).
\end{equation}
\item \begin{equation}
\label{R4}
 \mathbb{E}\left[ \left\| Z_t - Z_s \right\|^{2p} \right] \leq C_0 \left(1+ \mathbb{E}\left[\left\| Z_0 \right\|^{2p}\right]\right) \left(t - s\right)^{p}. 
\end{equation}
If $ \beta = 0$, we have a better estimation:
\begin{equation}
\label{R5}
\mathbb{E}\left[ \left\| Z_t - Z_s \right\|^{2p} \right] \leq C_0 \left(1+ \mathbb{E}\left[\left\| Z_0 \right\|^{2p}\right]\right) \left(t - s\right)^{2p}.
\end{equation}
\end{enumerate}
The constant $C_0$  only depends on $ \left\| \alpha(0) \right\|,  \left\| \beta(0) \right\|, \left\| \gamma(0)\right\|$, $T$, $p$, and the Lipschitz constants of $\alpha, \beta \text{and }\gamma$.
\end{aprop}
\begin{adem}
Only \eqref{R2} requires a proof, the other results are well known (see \cite{RY}). Let $t \in [0,h]$, and $s \in [0,t]$.  Applying a convexity inequality, we get:
\begin{equation*}
\begin{split}
 \left\| Z_s - Y_s \right\|^{2p} & = 3^{2p-1}\left(\left\| Z_0 - Y_0 \right\|^{2p} +  s^{2p-1}\displaystyle\int_{0}^{s} \left\|\alpha\left(Z_u\right) - \gamma \left(Y_u\right) \right\|^{2p} du  + \left\| \displaystyle\int_{0}^{s} \beta\left( Z_u \right) - \beta\left( Y_u \right)    dW_u  \right\|^{2p}\right).
\end{split}
\end{equation*}
Taking the expectation of the supremum, and using the Burkholder-Davis-Gundy inequality, we get:
\begin{equation*}
\begin{split}
\mathbb{E}\left[ \underset{s\leq t}{\sup}\left\| Z_s - Y_s \right\|^{2p} \right]  & \leq 3^{2p-1}\Bigg(\mathbb{E}\left[ \left\| Z_0 - Y_0 \right\|^{2p} \right] +  \left(2h\right)^{2p-1}\displaystyle\int_{0}^{t} \mathbb{E}\left[  \left\|\alpha\left(Z_u\right)  \right\|^{2p}\right]  + \mathbb{E}\left[  \left\|\gamma \left(Y_u\right) \right\|^{2p}\right] du  \\
& +Kh^{p-1} \displaystyle\int_{0}^{t} \mathbb{E}\left[\left\|  \beta\left( Z_u \right) - \beta\left( Y_u \right) \right\|^{2p} \right]   du  \Bigg)
\end{split}
\end{equation*} 
where $K$ is the constant that appears in the Burkholder-Davis-Gundy inequality. By the Lipschitz assumption, and by using \eqref{R1}, we obtain:
\begin{equation*}
\begin{split}
\mathbb{E}\left[ \underset{s\leq t}{\sup}\left\| Z_s - Y_s \right\|^{2p} \right]  & \leq 3^{2p-1}\Bigg(\mathbb{E}\left[ \left\| Z_0 - Y_0 \right\|^{2p} \right] +  R \left( \mathbb{E}\left[ 1+ \left\|Z_0  \right\|^{2p}\right]  + \mathbb{E}\left[ 1+ \left\|Y_0 \right\|^{2p}\right]  \right) h^{2p} \\
& +KT^{p-1} L^{2p} \displaystyle\int_{0}^{t} \mathbb{E}\left[ \underset{v\leq u}{\sup}\left\| Z_v - Y_v \right\|^{2p} \right]   du  \Bigg) 
\end{split}
\end{equation*} 
where $R =  4^{2p-1} \left(\max\left\{ L, \max\left\{\left| \alpha^i\left(0\right)\right|, i \in \left\{1,\ldots,n\right\}\right\},  \max\left\{\left| \gamma^i\left(0\right)\right|, i \in \left\{1,\ldots,n\right\}\right\} \right\}\right)^{2p} \exp\left(C_0T\right) $ and $L$ the common Lipschitz constant of $\alpha, \beta$, and $\gamma$. We conclude using Gronwall's lemma and changing $C_0$.
\end{adem} 

The next results are similar to Lemmas 2.5 and 2.6 in \cite{AJC}.  However, they require slightly different regularity assumptions since the Brownian vector fields commute.
\begin{alem}
\label{Moment_CC}
Assume that:
\begin{itemize}
\item  $\forall j \in \left\{0,\ldots,d\right\}, \sigma^j \in \mathcal{C}^1\left(\mathbb{R}^n,\mathbb{R}^n\right)$ with bounded first derivatives.
\item $\sum \limits_{j=1}^d \partial \sigma^j \sigma^j$ is a Lipschitz continuous function.
\end{itemize}
Then: $\forall p \ge 1,  \exists C^{\prime}_1 \in \mathbb{R}^*_+, \forall N \in \mathbb{N}^*, \forall t \in [0,T]$:
\begin{equation}
 \mathbb{E}\left[1 + \left\| \bar{X}^{0}_t\right\|^{2p} \right] \leq \exp(C^{\prime}_1 \check{\tau}_t ) \left(1+ \left\| x\right\|^{2p}\right), 
\end{equation}
\begin{equation}
\mathbb{E}\left[1 + \left\| \bar{X}_t\right\|^{2p} \right]  \leq \exp(C^{\prime}_1 \check{\tau}_t ) \left(1+ \left\| x\right\|^{2p}\right), 
\end{equation}
\begin{equation}
 \mathbb{E}\left[1 + \left\| \bar{X}^{d+1}_t\right\|^{2p} \right] \leq \exp(C^{\prime}_1 \check{\tau}_t ) \left(1+ \left\| x\right\|^{2p}\right). 
\end{equation}
\end{alem}
\begin{adem}
Let $p \ge 1$ and $ t \in [0,T]$, then $\exists k \in  \left\{0,\ldots,N-1\right\}$ such that $t_k < t \leq t_{k+1}$.
$\left(\bar{X}^{0}_s\right)_{t_k < s \leq t_{k+1}}$ is the solution of the following ODE:
\begin{equation*}
\left\{
    \begin{array}{ll}
 dZ_s  = \frac{1}{2} \sigma^0(Z_s) ds \\\\
  Z_{t_k} = X^{NV}_{t_{k}}.
\end{array}
\right.
\end{equation*}
$\left(\bar{X}_s\right)_{t_k < s \leq t_{k+1}}$ is the solution of the following SDE:
\begin{equation*}
\left\{
    \begin{array}{ll}
 dZ_s  = \frac{1}{2} \sum \limits_{j=1}^d \partial \sigma^j \sigma^j\left(Z_s\right) ds +  \sum \limits_{j=1}^d  \sigma^j(Z_s) dW^j_s \\
  Z_{t_k} = \bar{X}^{0}_{t_{k+1}}.
\end{array}
\right.
\end{equation*}
$\left(\bar{X}^{d+1}_s\right)_{t_k < s \leq t_{k+1}}$ is the solution of the following ODE:
\begin{equation*}
\left\{
    \begin{array}{ll}
 dZ_s  = \frac{1}{2} \sigma^0(Z_s) ds \\\\
  Z_{t_k} = \bar{X}_{t_{k+1}}.
\end{array}
\right.
\end{equation*}
Applying \eqref{R1} from the Proposition \ref{P1}, we get, $\forall t \in \left(t_k, t_{k+1} \right]$
\begin{equation}
 \mathbb{E}\left[ 1+ \left\| \bar{X}^{0}_t \right\|^{2p} \right] \leq \mathbb{E}\left[1+ \left\| X^{NV,\eta}_{t_{k}} \right\|^{2p}\right] \exp \left(C_0 h \right),
\label{b1} 
\end{equation}
\begin{equation}
 \mathbb{E}\left[ 1+ \left\| \bar{X}_t \right\|^{2p} \right] \leq  \mathbb{E}\left[1+\left\| \bar{X}^{0,\eta}_{t_{k+1}} \right\|^{2p}\right] \exp \left(C_0 h \right), 
\label{b2}
\end{equation}
\begin{equation}
 \mathbb{E}\left[ 1+  \left\| \bar{X}^{d+1}_t \right\|^{2p} \right] \leq \mathbb{E}\left[1+\left\| \bar{X}_{t_{k+1}} \right\|^{2p}\right] \exp \left( C_0 h \right). 
\label{b3}
\end{equation}
Using backward induction on \eqref{b1}, \eqref{b2} and \eqref{b3} we get:
\begin{equation}
 \mathbb{E}\left[ 1 + \left\| \bar{X}^{j}_t\right\|^{2p} \right] \leq  \exp\left(3 C_0 t_{k+1} \right) \left( 1 + \left\|x\right\|^{2p} \right) 
\end{equation}
for $j \in \left\{0,d+1\right\}$, and 
\begin{equation}
 \mathbb{E}\left[ 1 + \left\| \bar{X}_t\right\|^{2p} \right] \leq  \exp\left(3 C_0 t_{k+1} \right) \left( 1 + \left\|x\right\|^{2p} \right). 
\end{equation}
 \end{adem}
The proof of the next Lemma is a straightforward consequence of Lemma \ref{Moment_CC},  \eqref{R4} and \eqref{R5} from Proposition \ref{P1}. 
\begin{alem}
\label{Increment_CC}
Under the assumptions of the previous Lemma we have the following result.\\
$\forall p \ge 1,  \exists C^{\prime}_2 \in \mathbb{R}^*_+, \forall N \in \mathbb{N}^*, \forall t \in [0,T]$:
\begin{equation}
\label{_Bar_0_CC}
\mathbb{E}\left[ \left\| \bar{X}^{0}_t - X^{NV}_{\hat{\tau}_t}\right\|^{2p} \right] \leq   C^{\prime}_2 \left(1+ \left\| x\right\|^{2p}\right) h^{2p}, 
\end{equation}
\begin{equation}
 \mathbb{E}\left[ \left\| \bar{X}_t - \bar{X}^0_{\check{\tau}_t}\right\|^{2p} \right] \leq C^{\prime}_2 \left(1+ \left\| x\right\|^{2p}\right) h^{p}, 
\end{equation}
\begin{equation}
\mathbb{E}\left[ \left\| \bar{X}^{d+1}_t - \bar{X}_{\check{\tau}_t}\right\|^{2p} \right] \leq  C^{\prime}_2 \left(1+ \left\| x\right\|^{2p}\right) h^{2p}. 
\end{equation}
\end{alem}

\subsection{Proof of the strong convergence in the commutative case}

\begin{adem}
The error process is given by:
\begin{equation}
\begin{split}
X_s - X^{NV}_{s} & =  \frac{1}{2}  \displaystyle\int_{0}^{s}\sigma^0\left(X_u\right) -  \sigma^0 \left(\bar{X}^0_u\right)  du +   \frac{1}{2}  \displaystyle\int_{0}^{s}\sigma^0\left(X_u\right) -  \sigma^0 \left(\bar{X}^{d+1}_u\right) du  \\
& +  \sum \limits_{j=1}^d  \displaystyle\int_{0}^{s}  \left( \sigma^j\left(X_u\right) -  \sigma^j \left(\bar{X}_u\right) \right) \circ dW^j_u. 
\end{split}
\end{equation}
Then, using a convexity inequality, we get for $t \in [0,T]$ and $p\ge 1$
\begin{equation}
\label{Estim}
\mathbb{E}\left[ \underset{s\leq t}{\sup}\left\|X_t - X^{NV}_{t}\right\|^{2p} \right] \leq 3^{2p-1} \left(E_0 + E_{d+1} + E\right),
 \end{equation}
where
\begin{equation}
E _0=  \mathbb{E}\left[\underset{s\leq t}{\sup}\left\|  \frac{1}{2} \displaystyle\int_{0}^{s}  \sigma^0\left(X_u\right) -  \sigma^0 \left(\bar{X}^0_u\right)  du \right\|^{2p} \right],
\end{equation}
\begin{equation}
E_{d+1} =  \mathbb{E}\left[\underset{s\leq t}{\sup}\left\|  \frac{1}{2} \displaystyle\int_{0}^{s}  \sigma^0\left(X_u\right) -  \sigma^0 \left(\bar{X}^{d+1}_u\right)  du \right\|^{2p} \right],
\end{equation}
and
\begin{equation}
E =  \mathbb{E}\left[\underset{s\leq t}{\sup}\left\|  \displaystyle \sum \limits_{j=1}^d \int_{0}^{s} \left( \sigma^j\left(X_u\right) -  \sigma^j \left(\bar{X}_u\right) \right) \circ dW^j_u \right\|^{2p} \right].
\end{equation}
For the reader's convenience, the proof of this theorem is split into intermediate steps.

\textbf{Step 1: estimation of $E_0$}.\\\\
Introducing $\sigma^0 \left( X_{\hat{\tau}_u} \right)$ and $\sigma^0 \left(X^{NV}_{\hat{\tau}_u}\right)$, and using convexity inequality, we obtain:
\begin{equation*}
\begin{split}
E_0& \leq  \frac{3^{2p-1}}{2^{2p}}  \Bigg( \mathbb{E}\left[\underset{s\leq t}{\sup}\left\| \displaystyle\int_{0}^{s}  \sigma^0\left(X_u\right) -  \sigma^0 \left( X_{\hat{\tau}_u} \right)  du \right\|^{2p} \right]  
+ T^{2p-1} \displaystyle\int_{0}^{t} \mathbb{E}\left[  \left\| \sigma^0\left(X_{\hat{\tau}_u}\right) -  \sigma^0 \left(X^{NV}_{\hat{\tau}_u}\right) \right\|^{2p} \right] du\\
& + T^{2p-1} \displaystyle\int_{0}^{t} \mathbb{E}\left[ \left\|  \sigma^0 \left(X^{NV}_{\hat{\tau}_u}\right) -  \sigma^0 \left(\bar{X}^0_u\right) \right\|^{2p}  \right] du \Bigg).
\end{split}
\end{equation*}
The two last integrals are easy to estimate. On the one hand, since $\sigma^0$ is Lipschitz,
\begin{equation*}
 \displaystyle\int_{0}^{t} \mathbb{E}\left[  \left\| \sigma^0\left(X_{\hat{\tau}_u}\right) -  \sigma^0 \left(X^{NV}_{\hat{\tau}_u}\right) \right\|^{2p} \right] du \leq  L^{2p}  \displaystyle\int_{0}^{t}   \mathbb{E}\left[ \underset{v\leq u}{\sup}\left\|X_v - X^{NV}_{v}\right\|^{2p} \right] dv.
\end{equation*}
On the other hand, using \eqref{_Bar_0_CC} from Lemma \ref{Increment_CC} together with the Lipschitz property of $\sigma^0$, we get 
\begin{equation*}
\displaystyle\int_{0}^{t} \mathbb{E}\left[  \left\|  \sigma^0 \left(X^{NV}_{\hat{\tau}_u}\right) -  \sigma^0 \left(\bar{X}^0_u\right) \right\|^{2p} \right] du \leq L^{2p}  \displaystyle\int_{0}^{t} \mathbb{E}\left[   \left\|\bar{X}^0_u - X^{NV}_{\hat{\tau}_u} \right\|^{2p}\right]du\leq  L^{2p} C^{\prime}_2 T  \left(1 + \left\| x\right\| ^{2p}\right) h^{2p} .
\end{equation*}
Now we look at the first integral. For $i \in \left\{1,\ldots,n\right\}$, using the integration by parts formula, we have
\begin{equation*}
\displaystyle\int_{0}^{s}  \sigma^{i0} \left(X_u\right) - \sigma^{i0} \left( X_{\hat{\tau}_u} \right)  du = \displaystyle\int_{0}^{s}\left( \check{\tau}_u \wedge s - u\right)   d\left(\sigma^{i0} \left(X_u\right)\right).
\end{equation*}
Then, applying It\^o's formula, we get:
\begin{equation*}
\begin{split}
\displaystyle\int_{0}^{s}  \sigma^{i0}\left(X_u\right) - \sigma^{i0} \left( X_{\hat{\tau}_u} \right)  du &=  \displaystyle\int_{0}^{s} \left( \check{\tau}_u \wedge s - u\right) \left(\nabla \sigma^{i0}\left(X_u\right) . ~b\left(X_u\right) + \frac{1}{2}   tr\left(\sigma \left(X_u\right) \sigma^{*}\left(X_u\right) \nabla^2 \sigma^{i0} \left(X_u\right) \right) \right) du  \\
& + \displaystyle\int_{0}^{s} \left( \check{\tau}_u \wedge s - u\right) \sigma^{*}\left(X_u\right)\nabla \sigma^{i0}\left(X_u\right) .~dW_u  
\end{split}
\end{equation*}
where $\sigma = \left(\sigma^{ij}\right)_{1\leq i \leq n, 1\leq j\leq d}$ is the diffusion matrix. 
Taking the expectation of the supremum and using a convexity inequality, the Burkholder-Davis-Gundy inequality and $\check{\tau}_u - u \leq h$, we get:
\begin{equation*}
\begin{split}
\mathbb{E}\left[\underset{s\leq t}{\sup} \left| \displaystyle\int_{0}^{s}  \sigma^{i0}\left(X_u\right) - \sigma^{i0} \left( X_{\hat{\tau}_u} \right)  du \right |^{2p} \right] &\leq 3^{2p-1} h^{2p}\Bigg(  T^{2p-1} \displaystyle\int_{0}^{t}\mathbb{E}\left[ \left| \nabla \sigma^{i0}\left(X_u\right) .~ b\left(X_u\right)  \right|^{2p}  \right] du  \\
& +  \frac{ T^{2p-1}}{2^{2p}} \displaystyle\int_{0}^{t}\mathbb{E}\left[ \left|  tr\left(\sigma \left(X_u\right) \sigma^{*}\left(X_u\right) \nabla^2 \sigma^{i0} \left(X_u\right) \right) \right|^{2p}  \right] du \\
& + T^{p-1} K \displaystyle\int_{0}^{t}\mathbb{E} \left[ \left\| \sigma^{*}\left(X_u\right)\nabla \sigma^{i0}\left(X_u\right) \right\|^{2p} \right] du \Bigg)
\end{split}
\end{equation*}
where $K$ is the constant that appears in the Burkholder-Davis-Gundy inequality.
By the regularity assumption on $\sigma^j, j \in \left\{0,\ldots d\right\}$, we easily get two constants $\alpha_0\in \mathbb{R}_+$ and $q_1 \in \mathbb{N}^*$ which only depend on $p, T, \sigma$ and $\sigma^0$, such that
\begin{equation*}
\begin{split}
\mathbb{E}\left[\underset{s\leq t}{\sup} \left| \displaystyle\int_{0}^{s}  \sigma^{i0}\left(X_u\right) - \sigma^{i0} \left( X_{\hat{\tau}_u} \right)  du \right |^{2p} \right] &\leq \alpha_0 h^{2p} \displaystyle\int_{0}^{t}\mathbb{E}\left[1+  \left\| X_u\right\|^{2q_1} \right] du. 
\end{split}
\end{equation*}
Moreover, by the Lipschitz assumption on $b,\sigma^j, j \in \left\{0,\ldots d\right\}$, \eqref{R1} from Proposition \ref{P1} ensures that $\mathbb{E}\left[1+  \left\| X_u\right\|^{2q_1} \right] < + \infty$.
Finally, by combining our different inequalities, we obtain a constant $\beta_0 \in \mathbb{R}^*_+$ independent of $N$, such that 
\begin{equation}
\label{PC0}
E_0 \leq \beta_0   \left(h^{2p} + \displaystyle\int_{0}^{t}\mathbb{E}\left[ \underset{v\leq u}{\sup}\left\|X_v - X^{NV}_{v}\right\|^{2p} \right] du \right).
\end{equation}
\textbf{Step 2: estimation of $E_{d+1}$}.\\\\
Introducing $\sigma^0 \left( X_{\hat{\tau}_u} \right)$ and $\sigma^0 \left(X^{NV}_{\hat{\tau}_u}\right)$, and using convexity inequality, we get
\begin{equation*}
\begin{split}
E_{d+1} & \leq  \frac{3^{2p-1}}{2^{2p}}  \Bigg( \mathbb{E}\left[\underset{s\leq t}{\sup}\left\| \displaystyle\int_{0}^{s}  \sigma^0\left(X_u\right) -  \sigma^0 \left( X_{\hat{\tau}_u} \right)  du \right\|^{2p} \right] 
 +T^{2p-1} \displaystyle\int_{0}^{t} \mathbb{E}\left[  \left\| \sigma^{0}\left(X_{\hat{\tau}_u}\right) -  \sigma^0 \left(X^{NV}_{\hat{\tau}_u}\right) \right\|^{2p} \right] du \\
&+\mathbb{E}\left[\underset{s\leq t}{\sup}\left\| \displaystyle\int_{0}^{t}  \sigma^0\left(\bar{X}^{d+1}_u\right) -  \sigma^0 \left( X^{NV}_{\hat{\tau}_u} \right)  du \right\|^{2p} \right] \Bigg) .
\end{split}
\end{equation*}
We have already dealt with the two first terms int the right-hand side in the estimation of $E_{d+1}$.
Let $i \in \left\{1,\ldots,n\right\}$, we use integration by parts and It\^o's formula to get:
\begin{equation*}
\begin{split}
 \displaystyle\int_{0}^{s}  \sigma^{i0}\left(\bar{X}^{d+1}_u\right) -  \sigma^{i0} \left( X^{NV}_{\hat{\tau}_u} \right)   du = & \displaystyle\int_{0}^{s}\left(  s \wedge \check{\tau}_u - u\right)   d\left(\sigma^{i0} \left(\bar{X}^{d+1}_u\right)\right) + \displaystyle\int_{0}^{\check{\tau}_s} \left(  \check{\tau}_u - u \right) d\left(\sigma^{i0} \left(\bar{X}_u\right)\right) \\ &+ \displaystyle\int_{0}^{\check{\tau}_s} \left(  \check{\tau}_u - u \right) d\left(\sigma^{i0} \left(\bar{X}^0_u\right)\right)\\ 
= & \frac{1}{2}\displaystyle\int_{0}^{s}\left(  \check{\tau}_u - u \right)   \nabla \sigma^{i0} \left(\bar{X}^{d+1}_u\right).~\sigma^{0}\left(\bar{X}^{d+1}_u\right) du \\
&+ \sum \limits_{j=1}^d \displaystyle\int_{0}^{\check{\tau}_s} \left(  \check{\tau}_u - u \right) \left(\nabla \sigma^{i0} \left(\bar{X}_u\right) .~ \sigma^j\left(\bar{X}_u\right)\right) dW_u^j \\
&+ \frac12 \sum \limits_{j=1}^d \displaystyle\int_{0}^{\check{\tau}_s} \left(  \check{\tau}_u - u \right) \left(\nabla \sigma^{i0} \left(\bar{X}_u\right) .~ \partial \sigma^j \sigma^j\left(\bar{X}_u\right)\right) du\\
& + \frac12 \displaystyle\int_{0}^{\check{\tau}_s} \left(  \check{\tau}_u - u \right) tr\left(\sigma \sigma^* \nabla^2 \sigma^{i0} \left(\bar{X}_u\right)\right) du\\
 &+\frac{1}{2} \displaystyle\int_{0}^{\check{\tau}_s} \left(  \check{\tau}_u - u \right)\nabla\sigma^{i0} \left(\bar{X}^0_u\right).~ \sigma^0\left(\bar{X}^0_u\right) du. 
\end{split}
\end{equation*}
Once again, by the regularity assumption on $\sigma^j, \forall j \in \left\{1,\ldots,d\right\}$, and by using Lemma \ref{Moment_CC} and $\check{\tau}_u - u \leq h$, we get a constant $\alpha_2 \in \mathbb{R}_+$ independent of $N$ such that
\begin{equation*}
\mathbb{E}\left[\underset{s\leq t}{\sup} \left| \displaystyle\int_{0}^{s}  \sigma^{i0}\left(\bar{X}^{d+1}_u\right) -  \sigma^{i0} \left( X^{NV}_{\hat{\tau}_u} \right)   du \right |^{2p}  \right] \leq \alpha_2  h^{2p}.
\end{equation*}
Finally, by summing up our different inequalities, we obtain $\beta_{d+1} \in  \mathbb{R}_+^*$, independent of $N$, such that
\begin{equation}
\label{PC1}
E_{d+1}\leq  \beta_{d+1}  \left(h^{2p} + \displaystyle\int_{0}^{t}\mathbb{E}\left[ \underset{v\leq u}{\sup}\left\|X_v - X^{NV}_{v}\right\|^{2p} \right] du \right). 
\end{equation}
\textbf{Step 3: estimation of $E$}.\\\\
Re-expressing the integrals in It\^o's form, we get:
\begin{equation*}
\displaystyle \sum \limits_{j=1}^d \int_{0}^{s}  \left(\sigma^j\left(X_u\right) -  \sigma^j \left(\bar{X}_u\right) \right)\circ dW^j_u = \displaystyle \sum \limits_{j=1}^d  \int_{0}^{s} \left(\sigma^j\left(X_u\right) -  \sigma^j \left(\bar{X}_u\right)\right) dW^j_u  + \frac{1}{2} \displaystyle \sum \limits_{j=1}^d \int_{0}^{s}  \left(\partial \sigma^j \sigma^j\left(X_u\right) -  \partial \sigma^j \sigma^j \left(\bar{X}_u\right) \right)du.
\end{equation*}
Applying a convexity inequality, the Burkholder-Davis-Gundy inequality, together with the Lipschitz property of $\sigma^j, j \in \left\{1,\ldots,d\right\}$ and $\sum \limits_{j=1}^d \partial \sigma^j \sigma^j$,  we get 
\begin{equation*}
  \mathbb{E}\left[\underset{s\leq t}{\sup}\left\| \sum \limits_{j=1}^d \displaystyle\int_{0}^{s}   \left(\sigma^j\left(X_u\right) -  \sigma^j \left(\bar{X}_u\right) \right) \circ dW^j_u \right\|^{2p} \right] \leq  \alpha \displaystyle\int_{0}^{t} \mathbb{E}\left[
 \left\|X_u - \bar{X}_{u}\right\|^{2p} \right] du. 
\end{equation*}
where $\alpha = d\left(2L\right)^{2p}\left(KT^{p-1} + \frac{T^{2p-1}}{2^{2p}}\right) $ and $K$ is the constant that appears in the Burkholder-Davis-Gundy inequality. 
To estimate $\mathbb{E}\left[ \left\|X_u - \bar{X}_{u}\right\|^{2p} \right]$, we introduce the following piecewise continuous process $\left(\check{X}_u\right)_{u \in [0,T]}$ such that for $u \in [t_k,t_{k+1})$:
\begin{equation}
\check{X}_u = X_{t_k} +   \sum \limits_{j=1}^d \displaystyle\int_{t_k}^{u} \sigma^j \left(\check{X}_v\right) \circ dW^j_v.
\end{equation}
Then, using a convexity inequality, we obtain
\begin{equation*}
\mathbb{E}\left[\underset{s\leq t}{\sup}\left\|  \sum \limits_{j=1}^d  \displaystyle\int_{0}^{s} \sigma^j\left(X_u\right) -  \sigma^j \left(\bar{X}_u\right) \circ dW^j_u \right\|^{2p} \right] \leq  2^{2p-1}\alpha \displaystyle\int_{0}^{t} \mathbb{E}\left[
\left\|X_u - \check{X}_{u}\right\|^{2p} \right] + \mathbb{E}\left[ \left\|\check{X}_u - \bar{X}_{u}\right\|^{2p} \right] du. 
\end{equation*}
An estimation of the first expectation is given by \eqref{R2} in Proposition \ref{P1}:
\begin{equation*}
\mathbb{E}\left[ \left\|X_u - \check{X}_{u}\right\|^{2p} \right] \leq 2 C_0 \exp\left(C_0 h\right)  \mathbb{E}\left[ 1+ \left\| X_{t_k} \right\|^{2p} \right]h^{2p}.
\end{equation*}
We estimate the second expectation using \eqref{R3} in Proposition \ref{P1}:
\begin{equation*}
\mathbb{E}\left[ \left\|\check{X}_u - \bar{X}_{u}\right\|^{2p} \right] \leq \exp\left(C_0h\right) \mathbb{E}\left[ \left\|X_{t_k} - \bar{X}^0_{t_{k+1}}\right\|^{2p} \right].
\end{equation*}
Introducing $X^{NV}_{t_k}$ in the right-hand side of the inequality, we get:
\begin{equation*}
\mathbb{E}\left[ \left\|\check{X}_u - \bar{X}_{u}\right\|^{2p} \right] \leq 2^{2p-1} \exp\left(C_0T\right) \left(\mathbb{E}\left[ \left\|X_{t_k} - X^{NV}_{t_{k}}\right\|^{2p} \right] + \mathbb{E}\left[ \left\| \bar{X}^0_{t_{k+1}} - X^{NV}_{t_k} \right\|^{2p} \right] \right).
\end{equation*}
Applying \eqref{_Bar_0_CC} from Lemma \ref{Increment_CC} to the last expectation, we obtain:
\begin{equation*}
\mathbb{E}\left[ \left\|\check{X}_u - \bar{X}_{u}\right\|^{2p} \right] \leq 2^{2p-1} \exp\left(C_0T\right) \left(\mathbb{E}\left[ \underset{v\leq u}{\sup}\left\|X_v - X^{NV}_{v}\right\|^{2p} \right] + C^{\prime}_1\left(1 +  \left\| x \right\|^{2p} \right) h^{2p} \right).
\end{equation*}
We conclude by combining our different results, we get a constant $\beta \in \mathbb{R}_+^*$  independent of $N$, such that
\begin{equation}
\label{PC}
 E \leq \beta  \left(h^{2p} + \displaystyle\int_{0}^{t}\mathbb{E}\left[ \underset{v\leq u}{\sup}\left\|X_v - X^{NV}_{v}\right\|^{2p} \right] du \right). 
\end{equation}
\textbf{Step 4: conclusion}.\\\\
Combining \eqref{Estim}, \eqref{PC0}, \eqref{PC1},\eqref{PC}, and Gronwall's lemma, we get the following estimation.
\begin{equation*}
 \mathbb{E}\left[ \underset{t\leq T}{\sup}\left\|X_t - X^{NV}_{t}\right\|^{2p} \right] \leq 3^{2p-1}\left(\beta_0 + \beta_{d+1} + \beta \right) \exp\left( 3^{2p-1}\left(\beta_0 + \beta_{d+1} + \beta \right)T\right)  h^{2p}.
\end{equation*}
\end{adem}

\section{Appendix}
This section is devoted to the proof of Proposition \ref{ADNV-NV}. 
To compare $\hat{X}^{D,\eta}$ with $X^{NV,\eta}$, we introduce the following non-adapted interpolation of $\left(\hat{X}^{D,\eta}_{t_k}\right)_{0\leq k \leq N}$:
\begin{equation}
\label{_DNV}
\begin{split}
X^{D,\eta}_{t} & = \hat{X}^{D,\eta}_{\hat{\tau}_t} +  b\left(\hat{X}^{D,\eta}_{\hat{\tau}_t}\right) \Delta t +  \sum \limits_{j=1}^d \sigma^j\left(\hat{X}^{D,\eta}_{\hat{\tau}_t} \right) \Delta W^j_{t} + \frac{1}{2} \sum \limits_{j=1}^d \partial \sigma^j \sigma^j\left(\hat{X}^{D,\eta}_{\hat{\tau}_t}\right) \left( \left(\Delta W_{t}^j\right)^2 - \Delta t  \right)    \\
& +  \sum \limits_{\eta_t m < \eta_t j} \partial \sigma^j \sigma^m\left(\hat{X}^{D,\eta}_{\hat{\tau}_t}\right)\Delta W_{\check{\tau}_t}^m \Delta W_{t}^j.     
\end{split}
\end{equation}
Proposition \ref{ADNV-NV} is a consequence of the next lemma, which compares  $\hat{X}^{D,\eta}$ and $X^{D,\eta}$, and the next proposition, which compares $X^{D,\eta}$ and $X^{NV,\eta}$.
\begin{alem}
\label{DNVA}
Under the assumptions of Lemma \ref{DNVA_Lemma}, we have: 
\begin{equation}
\forall p \ge 1, \forall \epsilon>0, \exists C^{\prime}_D \in \mathbb{R}_+^*, \forall N \in \mathbb{N}^*, ~ \mathbb{E}\left[ \underset{t\leq T}{\sup}\left\|\hat{X}^{D,\eta}_{t} - X^{D,\eta}_{t}\right\|^{2p} \right] \leq C^{\prime}_D \frac{1}{N^{2p-\epsilon}}.
\end{equation}
\end{alem}
\begin{adem}
Let $p \ge 1, q,r > 1, \displaystyle \frac{1}{q} +\displaystyle\frac{1}{r} = 1$, and $t \in [0,T]$. Subtracting \eqref{_DNV} from \eqref{Approx_Discret}, we obtain
\begin{equation}
\begin{split}
\hat{X}^{D,\eta}_{t} - X^{D,\eta}_{t}  & = \sum \limits_{\eta_t m < \eta_t j} \partial \sigma^j \sigma^m\left(\hat{X}^{D,\eta}_{\hat{\tau}_t}\right)\Delta W_{t}^j \left( \Delta W_{t}^m  -  \Delta W_{\check{\tau}_t}^m\right) \\
& = - \sum \limits_{\eta_t m < \eta_t j} \partial \sigma^j \sigma^m\left(\hat{X}^{D,\eta}_{\hat{\tau}_t}\right)\Delta W_{t}^j \left( W_{\check{\tau}_t}^m -  W_{t}^m\right). 
\end{split}
\end{equation}
Then, combining a convexity inequality and the H\"{o}lder inequality, we get 
\begin{equation}
\begin{split}
\mathbb{E}\left[\underset{t\leq T}{\sup} \left\| \hat{X}^{D,\eta}_{t} - X^{D,\eta}_{t}\right\|^{2p} \right] & \leq \left(\frac{d^2-d}{2}\right)^{2p-1} \sum \limits_{ m < j}\mathbb{E}^{\frac{1}{q}}\left[  \underset{t\leq T}{\sup} \left\| \partial \sigma^j \sigma^m\left(\hat{X}^{D,\eta}_{\hat{\tau}_t}\right)\right\|^{2pq}\right]\mathbb{E}^{\frac{1}{r}}\left[  \underset{t\leq T}{\sup} \left | \Delta W_{t}^j \left( W_{\check{\tau}_t}^m -  W_{t}^m\right) \right |^{2pr}\right].
\end{split}
\end{equation}
Since $\partial \sigma^j \sigma^m$ for $j,m \in \left\{1,\ldots,d\right\}$, has an affine growth, using  Lemma \ref{DNVA_Lemma} we obtain a constant $\beta_3$, independent of $N$, such that:
\begin{equation}
\begin{split}
\mathbb{E}\left[\underset{t\leq T}{\sup} \left\| \hat{X}^{D,\eta}_{t} - X^{D,\eta}_{t}\right\|^{2p} \right] & \leq \beta_3 \sum \limits_{ m < j} \mathbb{E}^{\frac{1}{r}}\left[  \underset{t\leq T}{\sup} \left | \Delta W_{t}^j \left( W_{\check{\tau}_t}^m -  W_{t}^m\right) \right |^{2pr}\right].
\end{split}
\end{equation}
Then a straightforward calculation gives us:
\begin{equation}
\begin{split}
\mathbb{E}\left[\underset{t\leq T}{\sup} \left\| \hat{X}^{D,\eta}_{t} - X^{D,\eta}_{t}\right\|^{2p} \right] & \leq \beta_3  \sum \limits_{ m < j} \mathbb{E}^{\frac{1}{r}}\left[  \underset{k \in \left\{1,\ldots,N\right\} }{\sup} \underset{ ~ t_k < t \leq t_{k+1} }{\sup}   \left | \Delta W_{t}^j \left( W_{\check{\tau}_t}^m -  W_{t}^m\right) \right |^{2pr} \right]\\
& \leq \beta_3  \sum \limits_{ m < j}  \mathbb{E}^{\frac{1}{r}}\left[  \sum \limits_{k=1}^N \underset{ ~ t_k < t \leq t_{k+1} }{\sup}   \left | \left(W_{t}^j - W_{t_k}^j \right) \left( W_{t_{k+1}}^m -  W_{t}^m\right) \right |^{2pr} \right]\\
& \leq  \beta_3 N^{\frac{1}{r}}  \sum \limits_{ m < j}  \mathbb{E}^{\frac{1}{r}}\left[  \underset{ ~ 0 < t \leq t_{1} }{\sup}   \left | W_{t}^j \right|^{2pr} \underset{ ~ 0 < t \leq t_{1} }{\sup} \left| W_{t_1}^m - W_{t}^m \right |^{2pr} \right] \\
& \leq  \beta_3 N^{\frac{1}{r}}  \sum \limits_{ m < j}  \mathbb{E}^\frac{2}{r}\left[  \underset{ ~ 0 < t \leq t_{1} }{\sup}   \left | W_{t}^j \right|^{2pr}  \right] .
\end{split}
\end{equation}
Using Doob's submartingale inequality, we get
\begin{equation}
\begin{split}
\mathbb{E}\left[\underset{t\leq T}{\sup} \left\| \hat{X}^{D,\eta}_{t} - X^{D,\eta}_{t}\right\|^{2p} \right] & \leq \gamma_3 N^{\frac{1}{r}}   \sum \limits_{ m < j}  \mathbb{E}^\frac{2}{r}\left[ \left | W_{t_1}^j  \right|^{2pr}  \right]
\end{split}
\end{equation}
where $\gamma_3 = \beta_3 \left(\frac{r}{r-1}\right)^{2}$.
Finally, we obtain:
\begin{equation}
\begin{split}
\mathbb{E}\left[\underset{t\leq T}{\sup} \left\| \hat{X}^{D,\eta}_{t} - X^{D,\eta}_{t}\right\|^{2p} \right] & \leq C_D^{\prime} \frac{1}{N^{2p-\frac{1}{r}}}
\end{split}
\end{equation}
where $C_D^{\prime} = \displaystyle  \frac{d^2 -d}{2}\gamma_3 T^{2p} ~\mathbb{E}^\frac{2}{r}\left[\left|G\right|^{2pr}\right]$ and G a normal random variable.
\end{adem}

\begin{aprop}
\label{DNV}
Under the assumptions of Theorem \ref{EP_GC}, we have a first order of strong convergence:
\begin{equation}
\forall p \ge 1, \exists C^{\prime \prime}_D \in \mathbb{R}_+^*, \forall N \in \mathbb{N}^*, ~ \mathbb{E}\left[ \underset{t\leq T}{\sup}\left\|X^{NV,\eta}_{t} - X^{D,\eta}_{t}\right\|^{2p} \right] \leq C^{\prime \prime}_D h^{2p}.
\end{equation}
\end{aprop}

Before proving this proposition, we recall some useful results stated and proved in \cite{AJC}.
\begin{alem}
\label{Lemma0}
Assume that:
\begin{itemize}
\item $\forall j \in \left\{1,\ldots,d\right\},\sigma^j \in \mathcal{C}^{1}\left(\mathbb{R}^n,\mathbb{R}^n\right)$.
\item $\sigma^0, \sigma^j$ and $\partial \sigma^j \sigma^j, \forall j \in \left\{1,\ldots,d\right\}$, are Lipschitz continuous functions.
\item $F \in \mathcal{C}^2\left(\mathbb{R}^n,\mathbb{R}^n\right)$  with polynomially
growing first and second order derivatives.
\end{itemize}
Then, $\forall p \ge 1,  \exists C_{NV} \in \mathbb{R}^*_+, \forall N \in \mathbb{N}^*, $
\begin{enumerate}[(i)]
\item
\begin{equation}
\label{Moment_NV_S}
\forall t \in [0,T], \mathbb{E}\left[1 + \left\| X^{NV,\eta}_t\right\|^{2p} \right] \leq C_{1}  \left(1+ \left\| x\right\|^{2p}\right). 
\end{equation}
\item 
 \begin{equation}
\label{Moment_NV}
 \forall t \in [0,T], \forall j \in \left\{0,\ldots,d+1\right\}, \mathbb{E}\left[1 + \left\| \bar{X}^{j,\eta}_t\right\|^{2p} \right] \leq C_{1}  \left(1+ \left\| x\right\|^{2p}\right). 
\end{equation}
\item
 \begin{equation}
\label{Bar_j}
\forall t \in [0,T], \forall j \in \left\{0,\ldots,d+1\right\}, \mathbb{E}\left[ \left\| \bar{X}^{j,\eta}_t - X^{NV,\eta}_{\hat{\tau}_t}\right\|^{2p} \right] \leq C_{1} \left(1+ \left\| x\right\|^{2p}\right) h^{p}.
\end{equation}
\item 
\begin{equation}
\label{FS}
 \forall j \in \left\{0,\ldots,d+1\right\}, \mathbb{E}\left[ \underset{t\leq T}{\sup}\left\| \displaystyle\int_{0}^{t} F\left( \bar{X}^{j,\eta}_s\right) - F\left( X^{NV,\eta}_{\hat{\tau}_s}\right) ds  \right\|^{2p} \right]\leq C_1 h^{2p}. 
\end{equation}
\end{enumerate}
\end{alem} 

Now, we turn to the proof of  Proposition \ref{DNV}.

\begin{_adem} \textbf{of Proposition \ref{DNV} :}
Let $t \in [0,T]$ and $s \in [0,t]$.
Rewriting \eqref{_DNV} in integral form, we get
\begin{equation}
\label{DNV-Interpol_ITO}
\begin{split}
X^{D,\eta}_{t}& = x +  \displaystyle\int_{0}^{t} \sigma^0\left(X^{D,\eta}_{\hat{\tau}_s}\right) ds +  \frac{1}{2}\sum \limits_{j=1}^d  \displaystyle\int_{0}^{t} \partial \sigma^j \sigma^j\left(X^{D,\eta}_{\hat{\tau}_s}\right)  ds + \sum \limits_{j=1}^d  \displaystyle\int_{0}^{t}  \sigma^j\left(X^{D,\eta}_{\hat{\tau}_s}\right)  dW^j_s \\ 
& + \sum \limits_{j=1}^d  \displaystyle\int_{0}^{t} \partial \sigma^j \sigma^j\left(X^{D,\eta}_{\hat{\tau}_s}\right) \Delta W_{s}^j dW_s^j +  \displaystyle\int_{0}^{t} \sum \limits_{\eta_s m < \eta_s j} \partial \sigma^j \sigma^m\left(X^{D,\eta}_{\hat{\tau}_s}\right)\Delta W_{\check{\tau}_s}^m dW_s^j .
\end{split}
\end{equation}
Subtracting \eqref{DNV-Interpol_ITO} from \eqref{NV-Interpol_ITO}, we get
\begin{equation}
\begin{split}
X^{NV,\eta}_s - X^{D,\eta}_s  & =  \frac{1}{2} \displaystyle\int_{0}^{s} \sigma^0\left(\bar{X}^{0,\eta}_u\right) -  \sigma^0\left(X^{D,\eta}_{\hat{\tau}_u}\right) du + \sum \limits_{j=1}^d  \displaystyle\int_{0}^{s} \partial \sigma^j \sigma^j\left(\bar{X}^{j,\eta}_u\right) - \partial \sigma^j \sigma^j\left(X^{D,\eta}_{\hat{\tau}_u}\right)  du\\ 
& +  \frac{1}{2} \displaystyle\int_{0}^{s} \sigma^0\left(\bar{X}^{d+1,\eta}_u\right) -  \sigma^0\left(X^{D,\eta}_{\hat{\tau}_u}\right) du  + \sum \limits_{j=1}^d  \displaystyle\int_{0}^{s} \sigma^j\left(\bar{X}^{j,\eta}_u\right) -  \sigma^j\left(X^{D,\eta}_{\hat{\tau}_u}\right) dW_u^j \\
& - \sum \limits_{j=1}^d  \displaystyle\int_{0}^{s}  \left(\partial \sigma^j \sigma^j\left(X^{D,\eta}_{\hat{\tau}_u}\right) \Delta W_{u}^j + \sum \limits_{\eta_u m < \eta_u j} \partial \sigma^j \sigma^m\left(X^{D,\eta}_{\hat{\tau}_u}\right)\Delta W_{\check{\tau}_u}^m \right)    dW_u^j.
\end{split}
\end{equation} 
Using the Burkholder-Davis-Gundy inequality and a convexity inequality, we obtain
\begin{equation}
\label{Estimation}
\mathbb{E}\left[\underset{s\leq t}{\sup} \left\| X^{NV,\eta}_s - X^{D,\eta}_s \right\|^{2p} \right] \leq \left(2d+2\right)^{2p-1}\left(1+KT^{p-1}\right)\left( \sum \limits_{j=0}^{d+1} I_j + \sum \limits_{j=1}^{d} E_j \right)
\end{equation}
where $K$ is the constant that appears in the Burkholder-Davis-Gundy inequality,
\begin{equation*}
 E_j =\displaystyle\int_{0}^{t}\mathbb{E}\Bigg[ \Big\|\sigma^j\left(\bar{X}^{j,\eta}_u\right) - \sigma^j\left(X^{D,\eta}_{\hat{\tau}_u}\right) - \partial \sigma^j \sigma^j\left(X^{D,\eta}_{\hat{\tau}_u}\right) \Delta W_{u}^j  -  \sum \limits_{\eta_u m < \eta_u j} \partial \sigma^j \sigma^m\left(X^{D,\eta}_{\hat{\tau}_u}\right)\Delta W_{\check{\tau}_u}^m \Big\|^{2p} \Bigg] du 
\end{equation*}
and
\begin{equation*}
 I_j = \mathbb{E}\left[\underset{s\leq t}{\sup} \left\| \displaystyle\int_{0}^{s} F^j\left(\bar{X}^{j,\eta}_u\right) -  F^j\left(X^{D,\eta}_{\hat{\tau}_u}\right) du \right\|^{2p} \right] 
\end{equation*}
with $F^0 = F^{d+1} = \sigma^0$ and $F^j = \partial \sigma^j  \sigma^j, \forall j \in \left\{1,\ldots,d\right\}$.

\textbf{Step 1: estimation of $E_j$, for $j \in \left\{1,\ldots,d\right\}$}.\\\\
We first introduce the vector $R^j_u$, for $u \in [0,t]$, defined by:
\begin{equation*}
R^{j}_u =  \sigma^j\left(\bar{X}^{j,\eta}_u\right) - \sigma^j\left(X^{NV,\eta}_{\hat{\tau}_u}\right)  - \partial \sigma^j \sigma^j\left(X^{NV,\eta}_{\hat{\tau}_u}\right) \Delta W_{u}^j -  \sum \limits_{\eta_u m < \eta_u j} \partial \sigma^j \sigma^m\left(X^{NV,\eta}_{\hat{\tau}_u}\right)\Delta W_{\check{\tau}_u}^m, 
\end{equation*}
with coordinates $\left(R^{ij}_u\right)_{1\leq i \leq n}$.
Denoting 
\begin{equation}
J_u^{j,\eta} = \sum \limits_{\eta_u m <\eta_u j}  \sigma^m \left( X^{NV,\eta}_{\hat{\tau}_u}\right) \Delta W^m_{\check{\tau}_u} + \sigma^j \left( X^{NV,\eta}_{\hat{\tau}_u}\right) \Delta W^j_{u},
\end{equation}
$R^{ij}_u$ rewrites
\begin{equation}
\begin{split}
R^{ij}_u = \sigma^{ij}\left( \bar{X}^{j,\eta}_u\right)  - \sigma^{ij}\left( X^{NV,\eta}_{\hat{\tau}_u}\right) - \nabla \sigma^{ij}\left( X^{NV,\eta}_{\hat{\tau}_u}\right) .~ J^{j,\eta}_u.  
\end{split}
\end{equation}
By the mean value theorem, $ \exists \alpha_u^{ij} \in [0,1]$ such that:
\begin{equation}
\begin{split}
\sigma^{ij}\left( \bar{X}^{j,\eta}_u\right)  - \sigma^{ij}\left( X^{NV,\eta}_{\hat{\tau}_u}\right) = \nabla \sigma^{ij}\left(\xi^{ij}_u\right) .\left(\bar{X}^{j,\eta}_u -  X^{NV,\eta}_{\hat{\tau}_u}\right),
\end{split}
\end{equation}
with
\begin{equation}
\xi^{ij}_u  =  X^{NV,\eta}_{\hat{\tau}_u}  + \alpha_u^{ij} \left(\bar{X}^{j,\eta}_u -  X^{NV,\eta}_{\hat{\tau}_u}\right).
\end{equation}
Hence, introducing 
\begin{equation}
 \bar{R}^j_u =\bar{X}^{j,\eta}_u -  X^{NV,\eta}_{\hat{\tau}_u} -  J_u^{j,\eta},
\end{equation}
we get
\begin{equation}
R^{ij}_u =\nabla \sigma^{ij}\left( \xi^{ij}_u\right) .~ \bar{R}^j_u + \left(\nabla \sigma^{ij}\left( \xi^{ij}_u\right) -  \nabla \sigma^{ij}\left( X^{NV,\eta}_{\hat{\tau}_u}\right) \right). ~J_u^{j,\eta}.
\end{equation}
Let us now estimate $R^{ij}$. Since $\partial \sigma^j$ is Lipschitz continuous, we have 
\begin{equation}
\begin{split}
 \left| R^{ij}_u \right|^{2p} \leq  2^{2p-1} \left(1+L^{2p}\right) \left( \left\| \nabla\sigma^{ij}\left(\xi^{ij}_u\right)\right\|^{2p}\left\| \bar{R}^j_u\right\|^{2p} + \left\| \xi^{ij}_u -  X^{NV,\eta}_{\hat{\tau}_u}\right\|^{2p} \left\| J^{j,\eta}_u\right\|^{2p}   \right) 
\label{Rij}
\end{split}
\end{equation}
where $L \in \mathbb{R}_+^*$ denotes a common Lipschitz constant of $\sigma^j, j\in \left\{0,\ldots,d\right\}, \partial \sigma^j$ and $ \partial \sigma^j \sigma^m, j,m\in \left\{1,\ldots,d\right\}$. Then, using the Cauchy-Schwarz inequality, we get
\begin{equation}
\mathbb{E} \left[  \left\| \xi^{ij}_u -  X^{NV,\eta}_{\hat{\tau}_u}\right\|^{2p} \left\| J^{j,\eta}_u\right\|^{2p}\right] \leq \left(\mathbb{E} \left[  \left\| \xi^{ij}_u -  X^{NV,\eta}_{\hat{\tau}_u}\right\|^{4p}\right] \mathbb{E} \left[  \left\| J^{j,\eta}_u\right\|^{4p}\right]\right)^{\frac{1}{2}}.
\end{equation}
Applying \eqref{Bar_j} from Lemma \ref{Lemma0}, we obtain
\begin{equation}
\mathbb{E} \left[ \left\| \xi^{ij}_u -  X^{NV,\eta}_{\hat{\tau}_u}\right\|^{4p}\right] \leq C_{1} \left( 1 + \left\| x \right\|^{4p} \right) h^{2p}.
\end{equation}
Once again, combining the Cauchy-Schwarz inequality, the Lipschitz property for $\sigma^j, j \in \left\{1,\ldots,d\right\}$, and \eqref{Moment_NV_S} from Lemma \ref{Lemma0}, we get $\beta_1 \in \mathbb{R}_+$ independent of $N$ such that: 
\begin{equation}
\begin{split}
\mathbb{E} \left[ \left\| J^{j,\eta}_u\right\|^{4p} \right] & \leq \beta_1 h^{2p}.
\end{split}
\end{equation}
For the last term in the right-hand side of \eqref{Rij} we obtain:
\begin{equation}
\mathbb{E} \left[  \left\| \xi^{ij}_u -  X^{NV,\eta}_{\hat{\tau}_u}\right\|^{2p} \left\| J^{j,\eta}_u\right\|^{2p}\right] \leq \alpha_1 h^{2p}
\end{equation}
where $\alpha_1 = \left(\beta_1  C_{1} \left( 1 + \left\| x \right\|^{4p} \right)\right)^{\frac{1}{2}}$.
For the first term in the right-hand side of \eqref{Rij}, by the Lipschitz property of $\sigma^j$, $\partial \sigma^j $ is bounded by a constant denoted by $M$, so it remains to evaluate $\bar{R}^j$. Writing $\bar{R}^j$ in integral form, we obtain:
\begin{equation}
\begin{split}
 \bar{R}^j_u & =\bar{X}^j_u -  X^{NV,\eta}_{\hat{\tau}_u} -  J^{j,\eta}_u\\
& = \frac{1}{2} \displaystyle\int_{{\hat{\tau}_u}}^{\check{\tau}_u} \mathds{1}_{\left\{\eta_u = 1\right\}}\sigma^0\left(\bar{X}_v^0 \right) + \mathds{1}_{\left\{\eta_u = -1\right\}} \sigma^0\left(\bar{X}_v^{d+1} \right)dv +  \sum \limits_{\eta_u m < \eta_u j} \frac{1}{2} \displaystyle\int_{{\hat{\tau}_u}}^{\check{\tau}_u} \partial \sigma^m \sigma^m\left(\bar{X}_v^m \right) dv + \frac{1}{2} \displaystyle\int_{{\hat{\tau}_u}}^{u} \partial \sigma^j \sigma^j\left(\bar{X}_v^j \right) dv \\
& + \sum  \limits_{\eta_u m < \eta_u j}  \displaystyle\int_{{\hat{\tau}_u}}^{\check{\tau}_u} \sigma^m\left(\bar{X}_v^j \right) -  \sigma^m\left( X^{NV,\eta}_{\hat{\tau}_u} \right) dW^m_v +  \displaystyle\int_{{\hat{\tau}_u}}^{u} \sigma^j\left(\bar{X}_v^j \right) -  \sigma^j\left( X^{NV,\eta}_{\hat{\tau}_u} \right) dW^j_v.
\end{split}
\end{equation}
Combining a convexity inequality, the Burkholder-Davis-Gundy inequality and the Lipschitz property of $\sigma^m$ for $m \in \left\{0,\ldots,d\right\}$ and $\partial \sigma^m \sigma^m$ for $m \in \left\{1,\ldots,d\right\}$, together with \eqref{Moment_NV} and \eqref{Bar_j} from Lemma \ref{Lemma0}, we get $\beta_2 \in \mathbb{R}_+$ independent of $N$ such that: 
\begin{equation}
\begin{split}
\mathbb{E}\left[ \left\| \bar{R}^j_u  \right\|^{2p} \right]  & \leq \beta_2 h^{2p}. 
\end{split}
\end{equation}
For the first term in the right-hand side of \eqref{Rij}, we obtain:
\begin{equation}
\begin{split}
\mathbb{E}\left[ \left\| \nabla \sigma^{i,j}\left(\xi^{ij}_u\right)\right\|^{2p} \left\| \bar{R}^j_u  \right\|^{2p} \right]  & \leq \alpha_2 h^{2p} 
\end{split}
\end{equation}
where $\alpha_2 = M^{2p} \beta_2$. This leads us to the following estimation:
\begin{equation}
\begin{split}
\mathbb{E}\left[ \left| R^{ij}_u \right|^{2p} \right]  & \leq \alpha_0 h^{2p} 
\end{split}
\end{equation}
where $\alpha_0 = 2^{2p-1}\left(1+L^{2p}\right) \left(\alpha_1 + \alpha_2\right)$. Therefore
\begin{equation}
\mathbb{E}\left[ \left\| R^{j}_u \right\|^{2p} \right]  \leq n^{p} \alpha_0 h^{2p}. 
\end{equation}
Hence, by introducing $R^{j}_u $ in the expression of $E_j$ and by using  the Lipschitz assumption:
\begin{equation}
\begin{split}
& \mathbb{E}\Bigg[ \Big\|\sigma^j\left(\bar{X}^{j,\eta}_u\right) - \sigma^j\left(X^{D,\eta}_{\hat{\tau}_u}\right) -  \partial \sigma^j \sigma^j\left(X^{D,\eta}_{\hat{\tau}_u}\right) \Delta W_{u}^j - \sum \limits_{\eta_u m < \eta_u j} \partial \sigma^j \sigma^m\left(X^{D,\eta}_{\hat{\tau}_u}\right)\Delta W_{\check{\tau}_u}^m  \Big\|^{2p} \Bigg]\\
 & \leq  \left(d+2\right)^{2p-1} \left( 1+ L^{2p} \right) \Bigg( \mathbb{E}\left[ \underset{v\leq u}{\sup} \left\| X^{NV,\eta}_v - X^{D,\eta}_v \right\|^{2p} \right]+ \mathbb{E}\left[ \left\| X^{NV,\eta}_{\hat{\tau}_u} - X^{D,\eta}_{\hat{\tau}_u} \right\|^{2p} \left|\Delta W^j_u\right|^{2p}\right]\\ 
& +  \sum \limits_{m\neq j}  \mathbb{E}\left[ \left\| X^{NV,\eta}_{\hat{\tau}_u} - X^{D,\eta}_{\hat{\tau}_u} \right\|^{2p} \left|\Delta W^m_{\check{\tau}_u}\right|^{2p}\right] + \mathbb{E}\left[ \left\| R^{j}_u \right\|^{2p} \right] \Bigg).
\end{split}
\end{equation}
Then, by independence, for all $m \in \left\{1,\ldots,m\right\}$:
\begin{equation}
\begin{split}
 \mathbb{E}\left[ \left\| X^{NV,\eta}_{\hat{\tau}_u} - X^{D,\eta}_{\hat{\tau}_u} \right\|^{2p} \left|\Delta W^m_{\check{\tau}_u}\right|^{2p}\right] &  = \mathbb{E}\left[ \left\| X^{NV,\eta}_{\hat{\tau}_u} - X^{D,\eta}_{\hat{\tau}_u} \right\|^{2p}\right] \mathbb{E}\left[  \left|\Delta W^m_{\check{\tau}_u}\right|^{2p}\right]\\
 &\leq \mathbb{E}\left[\left|G\right|^{2p} \right] T^{p}~ \mathbb{E}\left[ \underset{v\leq u}{\sup} \left\| X^{NV,\eta}_v - X^{D,\eta}_v \right\|^{2p} \right]
\end{split}
\end{equation}
where $G$ is a normal random variable. 
Summing up these last inequalities, we get 
\begin{equation}
\label{Step1}
E_j \leq \gamma_1 \left(h^{2p}  + \displaystyle\int_{0}^{t}\mathbb{E} \left[\underset{v\leq u}{\sup} \left\| X^{NV,\eta}_{\hat{\tau}_v} - X^{D,\eta}_{\hat{\tau}_v} \right\|^{2p}\right] du\right) 
\end{equation}
where $\gamma_1 = \left(d+2\right)^{2p-1} \left( 1+ L^{2p} \right) \left(1 + d~\mathbb{E}\left[\left|G\right|^{2p} \right]T^p +n^p \alpha_0 T\right)$.

\textbf{Step 2: estimation of $I_j$, for $j \in \left\{0,\ldots,d+1 \right\}$}.\\\\
Let $j \in \left\{0,\ldots,d+1\right\}$.
\begin{equation*}
\begin{split}
\left\| \displaystyle\int_{0}^{s} F^j\left(\bar{X}^{j,\eta}_u\right) -  F^j\left(X^{D,\eta}_{\hat{\tau}_u}\right) du \right\|^{2p}  & \leq  2^{2p-1} \Bigg( \left\| \displaystyle\int_{0}^{s} F^j\left(\bar{X}^{j,\eta}_u\right) -  F^j\left(X^{NV,\eta}_{\hat{\tau}_u}\right) du \right\|^{2p}\\
 & + s^{2p-1} \displaystyle\int_{0}^{s}\left\| F^j\left(X^{NV,\eta}_{\hat{\tau}_u}\right) -  F^j\left(X^{D,\eta}_{\hat{\tau}_u}\right) \right\|^{2p} du  \Bigg).
\end{split}
\end{equation*}
Hence:
\begin{equation*}
I_j \leq \alpha_3 \left(\mathbb{E}\left[\underset{s\leq t}{\sup} \left\| \displaystyle\int_{0}^{s} F^j\left(\bar{X}^{j,\eta}_u\right) -  F^j\left(X^{NV,\eta}_{\hat{\tau}_u}\right) du \right\|^{2p} \right]  +\displaystyle\int_{0}^{t}\mathbb{E}\left[ \left\| F^j\left(X^{NV,\eta}_{\hat{\tau}_u}\right) -  F^j\left(X^{D,\eta}_{\hat{\tau}_u}\right) \right\|^{2p}\right] du \right)
\end{equation*}
where $\alpha_3 = 2^{2p-1}\left(1+  T^{2p-1}\right)$.
Then, by using \eqref{FS} from Lemma \ref{Lemma0} for the first integral and the Lipschitz assumption for the second one, we get:
\begin{equation}
\label{Step2}
I_j \leq \gamma_2 \left(h^{2p}  + \displaystyle\int_{0}^{t}\mathbb{E} \left[\underset{v\leq u}{\sup} \left\| X^{NV,\eta}_{\hat{\tau}_v} - X^{D,\eta}_{\hat{\tau}_v} \right\|^{2p}\right] du\right) 
\end{equation}
where $\gamma_2 =  \alpha_3 \left(C_1 + L^{2p}\right)$.

\textbf{Step 3: conclusion}\\\\
Finally, by combining \eqref{Step1}, \eqref{Step2}, together with \eqref{Estimation}, we obtain
\begin{equation}
\mathbb{E}\left[\underset{s\leq t}{\sup} \left\| X^{NV,\eta}_s - X^{D,\eta}_s \right\|^{2p} \right] \leq \gamma_3 \left(h^{2p}  + \displaystyle\int_{0}^{t}\mathbb{E} \left[\underset{v\leq u}{\sup} \left\| X^{NV,\eta}_{\hat{\tau}_v} - X^{D,\eta}_{\hat{\tau}_v} \right\|^{2p}\right] du\right)
\end{equation}
where $\gamma_3 = \left(2d+2\right)^{2p-1}\left(1+KT^{p-1}\right)\left(d \gamma_1 + \left(d+2\right)\gamma_2 \right)$ and we complete the proof using Gronwall's lemma since  Lemmas  \ref{DNVA_Lemma}, \ref{DNVA} and \ref{Lemma0} ensure that $\mathbb{E}\left[\underset{s\leq t}{\sup} \left\| X^{NV,\eta}_s - X^{D,\eta}_s \right\|^{2p} \right]$ is finite.
\end{_adem}

\end{document}